\documentclass[reqno,11pt]{amsart}
\usepackage{amssymb}
\usepackage{color}

\newtheorem{theorem}{Theorem}[section]
\newtheorem{lemma}[theorem]{Lemma}

\theoremstyle{definition}
\newtheorem{assumption}{Assumption}[section]

\theoremstyle{remark}
\newtheorem{remark}{Remark}[section]
\newcommand\bR{\mathbb{R}}
\newcommand\bQ{\mathbb{Q}}
\newcommand\cD{\mathcal{D}}
\newcommand\cF{\mathcal{F}}

\newcommand{\<}{\langle}
\renewcommand{\>}{\rangle}
\newcommand{\argmin}{\mathop{\sf argmin}}
\newcommand{\Trace}{\mathop{\sf trace}}
\numberwithin{equation}{section}
\begin{document}

\title[Diffusion approximation with
discontinuous coefficients] {On diffusion approximation with discontinuous
coefficients}

\author{N. V. Krylov}
\address{School of Mathematics, University of Minnesota,
Minneapolis, MN, 55455, USA}
\email{krylov@math.umn.edu}

\author{R. Liptser}
\address{Department of Electrical Engineering-Systems,
Tel Aviv University, 69978 Tel Aviv, Israel}
\email{liptser@eng.tau.ac.il}

\begin{abstract}
Convergence of stochastic processes with jumps to diffusion
processes is investigated in the case when the limit process has
discontinuous coefficients.
 An example is given
in which the diffusion approximation of a queueing
model yields a diffusion process with discontinuous
 diffusion and drift coefficients.
\end{abstract}



\maketitle

\section{Introduction}
\label{sec-1}

Suppose that we are given  a sequence of
semimartingales $(x^n_t)_{t\ge 0}$, $n=1,2,...$,
with paths in the Skorokhod space $\cD=
\cD([0,\infty),\bR^d)$ of $\bR^{d}$-valued right-continuous
functions on $[0,\infty)$
having left
limits on $(0,\infty)$. If one can prove that
the sequence of distributions $\bQ^{n}$ of $x^{n}_{\cdot}$
on $\cD$ weakly converges to the distribution
$\bQ$
of a diffusion process $(x _t)_{t\ge 0}$, then one says that
the sequence of $(x^n_t)_{t\ge 0}$ admits
a diffusion approximation. In this article by diffusion processes
we mean  solutions   of   It\^o equations
of the form
$$
x_t=x_0+\int_0^tb(s,x_s)\,ds+\int_0^t\sqrt{a(s,x_s)}\,dw_s,
$$
  with $w_{t}$ being
  a vector-valued Wiener process.
Usually to investigate the question if in a
particular situation there is a diffusion
approximation one uses
the  general framework of convergence
of semimartingales as developed for instance in \S3, Ch.~8 of
\cite{LS} (also see the references in this book).

The problem of diffusion approximation attracted
attention of
many researchers who obtained many deep and important results.
The reason for this is that
diffusion approximation is a quite
 efficient tool in stochastic systems
theory (see \cite{Ku'84}, \cite{Ku'90}), in
  asymptotic analysis of queueing models under heavy traffic and
bottleneck regimes  (see  \cite{KL}), in finding
asymptotically optimal filters (see
\cite{KuRu},  \cite{LR}), in asymptotical
optimization  in stochastic control problems (see
\cite{KuRu0},   \cite{LRT}), and in many other
issues.

In all above-mentioned references the coefficients
$a(t,x)$ and $b(t,x)$ of the limit diffusion process
are
continuous in $x$. In part, this is dictated by
the approach developed in \S3, Ch.~8 of \cite{LS}.
On the other hand, there are quite a few situations
in which the limit process should have discontinuous
coefficients. One of such situations is presented in
\cite{FS} where a queueing model is considered.
It was not possible to apply
standard results and the authors only
conjectured that the diffusion approximation
should be a process with natural coefficients.
Later this conjecture was rigorously proved in \cite{Ch}.
In \cite{Ch} and \cite{FS} only drift term is discontinuous.
Another
example of the limit diffusion  with discontinuous both
drift and diffusion
coefficients is given in article
\cite{KhasKryl} on averaging principle for diffusion processes
with null-recurrent fast component.

The idea to circumvent the discontinuity of $a$
and $b$ is to try to show that the time
spent by $(t,x_{t})$ in the set $G$ of their
discontinuity in $x$ is zero.
This turns out to be enough if
outside of $G$ the ``coefficients'' of $x^{n}_{t}$
converge ``uniformly'' to the coefficients of
$x_{t}$.
 By the way, even if all these hold,
still the functionals
$$
\int_{0}^{t}a(t,y_{t})\,dt,\quad
\int_{0}^{t}b(t,y_{t})\,dt,\quad y_{\cdot}\in\cD
$$
need not be continuous on the support of $\bQ$.
This closes the route of ``trivial''
generalizing the result from
\S3, Ch.~8 of \cite{LS}.

To estimate the time spent by $x_{t}$ we use
an inequality similar to the following
one
\begin{equation}
                                          \label{*}
E\int_{0}^{T}f(t,x_{t})\,dt \leq
N\Bigg(\int_{0}^{T}\int_{\mathbb{R}^d}
f^{d+1}(t,x)\,dxdt\Bigg)^{1/(d+1)},
\end{equation}
which is obtained in \cite{Kr74} for nonnegative
Borel $f$.
Then upon assuming that
$G\subset(0,\infty)\times\bR^{d}$
has $d+1$-dimensional Lebesgue measure zero
and substituting $I_{G}$ in place of $f$
in (\ref{*}) we get that indeed the time spent by $(t,x_{t})$
in $G$ is zero.
However, for (\ref{*}) to hold we need the process $x_{t}$
to be uniformly nondegenerate which may be not convenient
in some applications. Therefore,
in Sec.~\ref{section 3.14.1} we prove
a version of (\ref{*}), which
allows us to get
the conclusion about the time spent in $G$
assuming that the process is nondegenerate only on $G$.
In essence, our approach to diffusion approximation
with discontinuous coefficients is close to the one from
\cite{Ch}. However, details are quite different and
we get more general results
under less restrictive assumptions.
In particular, we do not impose
the linear growth condition. Neither do we assume that
the second moments of $x^{n}_{0}$ are bounded.
The weak limits
of processes with jumps appear in many other
settings, in particular, in Markov chain
approximations in the theory of controlled diffusion processes,
where, generally, the coefficients of $x^{n}_{t}$
are not supposed to converge to anything in any sense and yet
the processes converge weakly to a process of diffusion
type.

We mention here Theorem 5.3 in Ch.~10 of \cite{KD}
also bears on this matter
in the particular case of Markov chain
approximations in the theory of controlled diffusion processes.
Clearly, there is no way to specify
precisely the coefficients of all limit points in the general problem.
 Still one can obtain
some nontrivial information
and one may wonder if one can get anything from general results
when we are additionally
given that the coefficients do converge
on the major part of the space.
In Remarks \ref{remark 10.13.1} and \ref{remark 10.13.2}
we show that this is not the case in what concerns
Theorem 5.3 in Ch.~10 of \cite{KD}.

Above we alluded to the ``coefficients'' of $x^{n}_{t}$.
By them we actually mean the local drift and
the matrix of quadratic variation. We do not use
any additional structure of $x^{n}_{t}$. In particular,
the quadratic variation is just the sum of two terms:
one coming  from diffusion and another from
jumps. Therefore unlike \cite{KP} we do not use
any stochastic equations for $x^{n}_{t}$. This allows us
to neither  introduce nor use
  any assumptions on the martingales
driving these equations and their (usual) coefficients
thus making the presentation simpler and more general.
On the other hand it is worth noting that
the methods of \cite{KP} may be more useful
in other problems. Our intention was not to cover
all aspects of diffusion approximation
but rather give a new method allowing us to treat
discontinuous coefficients. In particular,
we do not discuss uniqueness of solutions to the limit
equation. This is a separate issue belonging to the theory
of diffusion processes and we only mention article \cite{KhasKryl},
where the reader can find a discussion of it.

The paper is organizes as follows. In Section \ref{section 4.14.1}
we prove our main results,
Theorems \ref{theorem 3.8.1} and \ref{theorem 4.25.1},
 about diffusion approximation. Their proofs
rely on the estimate proved in Sec.~\ref{section 3.14.1}
we have been talking about above.
But even if the set $G$ is empty,
the results which we prove are the first ones
of the kind.

In Theorems \ref{theorem 3.8.1} and \ref{theorem 4.25.1}
there is no assumption about any control
  of $\sqrt{a(t,x)}$ and $b(t,x)$
as $|x|\to\infty$, but instead we assume that
$\bQ^{n}$ converge weakly to $\bQ$. Therefore,
in Sec.~\ref{section 4.14.2}
we give a sufficient condition
for precompactness of a sequence of
distributions on Skorokhod space.
Interestingly enough, this condition
is  different from those which one gets from
\cite{JS} and \cite{LS} and again
does not involve usual growth conditions.
 Sec.~\ref{section 4.18.1} contains an example
of application of our results to
a queueing model close to the one from \cite{Ch}, \cite{FS}.
We slightly modify the model from \cite{Ch}, \cite{FS}
and get the diffusion approximation  with
discontinuous {\em drift} and {\em diffusion\/} coefficients.
To the best of our knowledge this is the first
example when the diffusion approximation leads
to discontinuous diffusion coefficients.

  The authors are sincerely grateful to the referees
for many useful suggestions.

\section{The main results}
                                            \label{section 4.14.1}
We use notions and notation from \cite{LS}.
For each $n=1,2,...$, let
$$
(\Omega^{n},\cF^{n},\cF^{n}_{t},t\geq0, P^{n})
$$
be a  stochastic basis
satisfying the ``usual'' assumptions.
Let $\cD$ be the Skorokhod space or right-continuous
$\bR^{d}$-valued functions $x_{t}$ given on $[0,\infty)$
and having left limits on $(0,\infty)$. As usual we endow
$\cD$ with Skorokhod-Lindvall metric in which
$\cD$ becomes a Polish space (see Theorem 2,  \S1, Ch.~6
of \cite{LS}).

Suppose that for each $n$
on $\Omega^{n}$ we are given an
$\cF^{n}_{t}$-semimartingale $x^{n}_{t}$, $t\geq0$,
with trajectories in $\cD$.
Let $(B^{n},C^{n},\nu^{n})$ be the triple of
predictable characteristics of $(x^{n}_{t},\cF^{n}_{t})$
and $\mu^{n}$ be its jump measure
(see \S1, Ch.~4 of \cite{LS}).
  Then
$$
x^{n}_{t}=x^{n}_{0}+B^{n}_{t}+x^{nc}_{t}+\int_{0}^{t}
\int_{|x|\leq1}x\,(\mu^{n}-\nu^{n})(dsdx)+\int_{0}^{t}\int_{|x|>1}
x\,\mu^{n}(dsdx),
$$
where
$B^{n}_{t}$ is a predictable process of locally bounded
variation with $B^{n}_{0}=0$, $x^{nc}_{t}$ is a continuous
local martingale with $\< x^{nc}\>_{t}=C^{n}_{t}$,
$\nu^{n}$ is the compensator of $\mu^{n}$.
Define
$$
m^{n }_{t}=x^{nc}_{t}+\int_{0}^{t}\int_{|x|\leq1}
x\,(\mu^{n}-\nu^{n})(dsdx),\quad
j^{n}_{t}=\int_{0}^{t}\int_{|x|>1}
x\,\mu^{n}(dsdx)
$$
so that $m^{n}_{t}$ is a locally square-integrable martingale
and
\begin{equation}\label{4.24.2}
x^{n}_{t}=x^{n}_{0}+B^{n}_{t}+m^{n}_{t}+j^{n}_{t}.
\end{equation}

\begin{assumption}{\rm
                                         \label{assumption 3.8.2}
(i) For each
$n$ on $(0,\infty)\times\cD$ we are given an
$\bR^{d}$-valued function
$b^{n}=b^{n}(t,y_{\cdot})$ and a $d\times d$
matrix valued function
$a^{n}=a^{n}(t,y_{\cdot})$ which is nonnegative
and symmetric for any $t$
and $y_{\cdot}\in\cD$. The functions $b^{n}$ and $a^{n}$
 are Borel measurable.
(ii) For each $r\in[0,\infty)$ there exists a locally integrable
function $L(r,t)$ given on $[0,\infty)$ such that
$L(r,t)$ increases in $r$ and
\begin{equation}
                                               \label{4.14.3}
|b^{n}(t,y_{\cdot})|+\Trace\,a^{n}(t,y_{\cdot})\leq L(r,t)
\end{equation}
whenever $t>0$, $y_{\cdot}\in\cD$, and $|y_{t}|\leq r$.
(iii) We have
$$
B^{n}_{t}=\int_{0}^{t}b^{n}(s,x^{n}_{\cdot})\,ds,\quad
\< m^{n}\>_{t}=2\int_{0}^{t}a^{n}(s,x^{n}_{\cdot})\,ds.
$$
}\end{assumption}

\begin{remark}{\rm
We have
$$
\< m^{n}\>^{ij}_{t}=
\< x^{nc}\>^{ij}_{t}+\int_{0}^{t}
\int_{|x|\leq1}x^{i}x^{j} \nu^{n}(dsdx)
$$
and it follows from Assumption \ref{assumption 3.8.2}
that both summands on the right are absolutely continuous
in $t$.
In particular, they are continuous, which along
with the continuity of $B^{n}_{t}$ implies that
$x^{n}_{t}$ is quasi leftcontinuous (see
Theorem 1, \S1, Ch.~4 of \cite{LS}).
}\end{remark}

\begin{assumption}{\rm
                                    \label{assumption 3.8.3}
(i)
On $(0,\infty)\times\bR^{d}$ we are given
an $\bR^{d}$-valued function $b=b(t,x)$ and
a $d\times d$ matrix valued function $a=a(t,x)$
which is nonnegative and symmetric for any $t$ and $x$.
The functions $b$ and $a$ are Borel measurable.

(ii) There exists a Borel
 set $G\subset(0,\infty)\times\bR^{d}$
(perhaps empty) such that, for almost
every $t\in(0,\infty)$, for every $x$
lying outside of the $t$-section
$G_{t}:=\{x\in\bR^{d}:(t,x)\in G\}$
of $G$ and any   sequence  $y^{n}_{\cdot}\in\cD$,
which converges to a continuous
function $y_{\cdot}$ satisfying
  $y _{t}=x$, it holds that
$$
b^{n}(t,y^{n}_{\cdot})\to b(t,x),\quad
a^{n}(t,y^{n}_{\cdot})\to a(t,x).
$$
}\end{assumption}
\begin{remark}{\rm
                                        \label{remark 3.13.2}
It is easy to see that Assumption \ref{assumption 3.8.3}
implies that for almost any $t$, the functions
$a(t,x)$ and $b(t,x)$ are continuous on the set
$\bR^{d}\setminus G_{t}$ in the relative topology of this set.

Also,
  Assumptions \ref{assumption 3.8.2} and \ref{assumption 3.8.3}
obviously imply that
$$
|b (t,x)|+\Trace\,a (t,x)\leq L(r,t)
$$
for almost
every $t\in(0,\infty)$ and all $x$ satisfying
$ |x|\leq r$, $x\not\in G_{t}$.
}\end{remark}
\begin{assumption}{\rm
                                         \label{assumption 3.8.4}
If $G\ne\emptyset$, then for almost each $t$

(i)  the set $ G _{t} $ has
Lebesgue measure zero,

(ii) for every $x\in  G _{t} $
and each sequence  $y^{n}_{\cdot}\in\cD$,
which converges to a continuous function $y_{\cdot}$ satisfying
  $  y _{t} =x$,
we have
\begin{equation}
                                              \label{4.19.2}
\varliminf_{n\to\infty}\det a^{n}(t,y^{n}_{\cdot})
\geq\delta(t,x)>0,
\end{equation}
where $\delta$ is a Borel function.
}\end{assumption}

\begin{remark}{\rm
Condition (\ref{4.19.2}) is satisfied if, for instance,
the processes $x^{n}_{t}$ are uniformly nondegenerate in
a neighborhood of $G_{t}$.
}\end{remark}

\begin{assumption}{\rm
                                         \label{assumption 3.8.5}
For any $T,\varepsilon\in(0,\infty)$, and any
$\alpha\in(0,1]$, it holds that
$$
\lim_{n\to\infty}P^{n}\big(
\nu^{n}\big((0,T]\times B^{c}_{\alpha}))\geq\varepsilon\big)=0,
$$
where $B_{\alpha}=\{x\in\bR^{d}:|x|<\alpha\}$,
 $B^{c}_{\alpha}=\{x\in\bR^{d}:|x|\geq \alpha\}$.
}\end{assumption}
\begin{remark}{\rm
                                        \label{remark 3.13.1}
Notice that for each $\alpha\in(0,1]$
and $r,T\in[0,\infty)$

$$
\theta^{n}_{rT}:=\int_{0}^{T}
\int_{|x|\leq1}|x|^{3}I_{|x_{s}|\leq r }
\,\nu^{n}(dsdx)\leq\int_{0}^{T}\int_{|x|<\alpha}+
\int_{0}^{T}\int_{|x|\geq\alpha}
$$
$$
\leq\alpha\int_{0}^{T}\int_{|x|\leq1}|x|^{2}I_{|x_{s}|\leq r }
\,\nu^{n}(dsdx)+\nu^{n}\big((0,T]\times B^{c}_{\alpha})),
$$
where according to Assumption \ref{assumption 3.8.2}
the first term on the right is less than
$$
 2\alpha\int_{0}^{T}I_{|x_{s}|\leq r }
\Trace\,a^{n}(s,x^{n}_{\cdot})
\,ds\leq 2\alpha\int_{0}^{T}L(r,s)\,ds.
$$
It follows easily that, owing to
Assumptions  \ref{assumption 3.8.2}
and \ref{assumption 3.8.5},
for each $\varepsilon>0$ and $r,T\in[0,\infty)$, we have
$$
\lim_{n\to\infty}
P^{n}(\theta^{n}_{rT}\geq\varepsilon\big)=0
$$
and since $\theta^{n}_{rT}\leq2\int_{0}^{T}L(r,s)\,ds$,
we also have $E^{n}\theta^{n}_{rT}\to0$ as $n\to\infty$,
where $E^{n}$ is the expectation sign relative to $P^{n}$.
}\end{remark}
\begin{remark}{\rm
                                         \label{remark 4.17.1}
Define
\begin{equation}
                                            \label{4.17.2}
\gamma^{n}=\inf\{t\geq0:|j^{n}_{t}|>1\}.
\end{equation}
Then $\gamma^{n}$ is an $\cF^{n}_{t}$-stopping time,
and obviously $j^{n}_{t}=0$ for $0\leq t<\gamma^{n}$.
Furthermore,
by Lemma VI.4.22 of \cite{JS},
Assumption \ref{assumption 3.8.5} implies that
$$
P^{n}(\gamma^{n}\leq T)\to0
$$
for each $T\in[0,\infty)$.
}\end{remark}
\begin{theorem}
                                          \label{theorem 3.8.1}
In addition to Assumptions
\ref{assumption 3.8.2}-\ref{assumption 3.8.5},
 suppose that
 the sequence of distributions $(\bQ^{n})_{n\geq1}$ of
$x^{n}_{\cdot}$
converges weakly on the Polish space
 $\cD$  to a measure
$\bQ$. Then $\bQ$ is the distribution
 of a solution  of the It\^o equation
\begin{equation}
                                                  \label{4.18.3}
x_{t}=x_{0}+\int_{0}^{t}\sqrt{2a(s,x_{s})}\,dw_{s}
+\int_{0}^{t}b(s,x_{s})\,ds
\end{equation}
defined on a probability space with $w_{t}$
being a $d$-dimensional
Wiener process.
\end{theorem}

\begin{remark}{\rm
                                             \label{remark 10.13.1}
Notice that there are {\em no conditions\/} on the values of
$ a(t,x)$ and $b(t,x)$ on the set $G$. Hence Theorem
\ref{theorem 3.8.1} holds if we replace
 $ a,b  $ with any other Borel functions,
which coincide with the original ones on the complement
$\Gamma$ of $G$.
Of course, this can only happen if
$$
\int_{0}^{t}I_{G}(s,x_{s})\, ds=0\quad\text{(a.s.)}.
$$
This equality is proved in Lemma \ref{lemma 4.21.4}. In particular,
$x_{t}$ satisfies
\begin{equation}
                                                  \label{4.18.30}
x_{t}=x_{0}+\int_{0}^{t}I_{\Gamma}(s,x_{s})\sqrt{2a(s,x_{s})}\,dw_{s}
+\int_{0}^{t}I_{\Gamma}(s,x_{s})b(s,x_{s})\,ds.
\end{equation}
Thus, the limit process satisfies (\ref{4.18.30}).
A particular feature of this equation is
that generally its solutions   are not  unique.
Indeed, let $x'_{t}$ be a one-dimensional Wiener process $w_{t}$ and
$x''_{t}$ the process identically equal to zero. They both satisfy
$dx_{t}=\sqrt{2a(t,x_{t})}\,dw_{t}$, where $a(t,x)=1/2$ for
$(t,x)\not\in G$,
$a(t,x)=0$ for $(t,x)\in G$, and $G=[0,\infty)\times\{0\}$. Of course,
there are many more different solutions which spend some time at zero
then follow the trajectories of $w_{t}$ for a while and then again
stay at zero. Therefore, the statement that $x_{t}$
has the form
$$
x_{t}=x_{0}+\int_{0}^{t}\sqrt{2a_{s}}\,dw_{s}
+\int_{0}^{t}b_{s}\,ds,
$$
where $a_{s}=a(s,x_{s})$ and $b_{s}=b(s,x_{s})$
whenever $(s,x_{s})\not\in G$ and $a$ and $b$ are not specified
otherwise (cf. the first part of
Theorem 5.3 in Ch.~10 of \cite{KD}), contains
very little information on the process: in the above example
both $x'_{t}$ and $x''_{t}$ have this form.
In contrast with this always in the above example, the fact that
without changing $x_{t}$ one can change $a,b$ on $G$ in any way,
and thus take $a\equiv1/2$, leaves
only one possibility: $x_{t}=w_{t}$.

}\end{remark}
\begin{remark}{\rm
                                             \label{remark 10.13.2}
From Remark \ref{remark 10.13.1}
we also see that the assumption that (\ref{4.18.3}) has a unique
(weak or strong) solution makes no sense unless
the values of $ a(t,x)$ and $b(t,x)$ are {\em specified everywhere\/}.
In Theorem 5.3 in Ch.~10 of \cite{KD} an attempt  is presented
to specify $ a(t,x)$ and $b(t,x)$ on $G$ consisting of
requiring that they belong to the set of all possible
diffusion and drift coefficients of $x_{t}$ when $x_{t}\in
G_{t}$. Generally, the   set  $x_{t}\in
G_{t}$ has zero probability
 (say, for the Wiener process) and  the requirement
seems to have little sense. Nevertheless, it is natural
to assume that, if $x_{t}=w_{t}$ in the example from
Remark \ref{remark 10.13.1}, then the only possibility for
$a(t,0)$ is $1/2$, the same value as for all other $x$.

In that case, the equation $dx_{t}=\sqrt{2a(t,x_{t})}\,dw_{t}$
($=dw_{t}$)
with zero initial condition has a unique solution,
the distribution of which (by Theorem \ref{theorem 3.8.1}) is the weak
limit of the distributions of solutions to $dx^{n}_{t}=
\sqrt{2a^{n}(x_{t})}\,dw_{t}$ with zero initial condition,
where $a^{n}(x)=1/2$ for $|x|\geq1/n$ and $a^{n}(x)=1/3$
for $|x|<1/n$.

Hovewer, this fact does not imply that
the distributions of
  any other
processes $z^{n}_{\cdot}$ converge to the Wiener measure,
provided only that $z^{n}_{t}$
satisfy $z^{n}_{0}=0$ and
 $dz^{n}_{t}=\sqrt{2c^{n}(z^{n}_{t})}\,dw_{t}$
with $c^{n}(x)=a^{n}(x)$ for $|x|\geq1/n$, $c^{n}\geq0$, and
$\sup_{n,x}c^{n}(x)<\infty$.
To show this, it suffices to define $c^{n}(x)=n^{2}x^{2}$
for $|x|\leq1/n$ and notice that $z^{n}_{t}\equiv0$ for all $n$.

This somewhat contradicts the second part of
Theorem 5.3 in Ch.~10 of \cite{KD}.

}\end{remark}

The proof of Theorem \ref{theorem 3.8.1}
consists of several steps  throughout which we assume
that the conditions of this theorem are satisfied.

The idea
is to rewrite  (\ref{4.18.3})  in terms of the
martingale problem of Stroock-Varadhan.
Then naturally we also want to write the information
about $x^{n}_{t}$ in a
martingale form not involving stochastic
bases and convenient to passing to the limit.
This is done in  Lemma
\ref{lemma 4.21.1}.
After that we pass to the limit and
in Lemma \ref{lemma 4.21.2} derive our theorem
upon additionally assuming that the
time spent by the limit process
$(t,x_{t})$ in the set $G$ of possible discontinuities
of coefficients is zero. This additional assumption holds,
for instance, if $G=\emptyset$.
 Lemma \ref{lemma 4.21.4} concludes the proof
of the theorem.

After that in Theorem \ref{theorem 4.25.1} we extend
Theorem \ref{theorem 3.8.1}
 to cases in which uniform nondegeneracy on $G_{t}$
of diffusion is not required.
We show the usefulness of Theorem \ref{theorem 4.25.1}
in Remark \ref{remark 4.25.1}.

As any probability measure on $\cD$, the measure
$\bQ$ is the distribution on $\cD$ of a process $x_{\cdot}$
having trajectories in $\cD$ and defined on a probability
space. By $E$ we denote the expectation sign
associated with that probability
space. We will see that the theorem holds for this
$x_{\cdot}$ up to a possible enlargement of the probability
space on which $x_{\cdot}$ lives.
In the following lemma Assumptions \ref{assumption 3.8.3} and
\ref{assumption 3.8.4} are not used.

By $C^{\infty}_{0}(\bR^{d+1})$
we denote the set of all
infinitely differentiable real-valued function $u=u(t,x)$
on $\bR^{d+1}$ with compact support.

\begin{lemma}
                                        \label{lemma 4.21.1}
For any $0\leq t_{1}\leq...\leq t_{q}\leq s\leq
t<\infty$,
 continuous bounded  function $f$ on $\bR^{qd}$,
and $u\in C^{\infty}_{0}(\bR^{d+1})$, we have
\begin{eqnarray}
&Ef(x _{t_{1}},...,x _{t_{q}})
\big[u(t,x _{t})-u(s,x _{s})\big]&\nonumber
\\
&=\lim_{n\to\infty}E^{n}f(x^{n}_{t_{1}},...
,x^{n}_{t_{q}})
\int_{s}^{t}\big[u_{p}(p,x^{n}_{p})&+a^{nij}(p,x^{n}_{\cdot})
u_{x^{i}x^{j}}(p, x^{n}_{p})\nonumber
\\
                                              \label{4.21.6}
& &+b^{ni}(p,x^{n}_{\cdot})u_{x^{i}}(p,x^{n}_{p})\big]\,dp.
\end{eqnarray}
Furthermore, the integrand with respect to $p$ is less
than $NL(r,p)$, where the constants
 $N$ and $r$ depend only on $u$
but not on $\omega$ and $n$.
\end{lemma}

Proof. Denote
$$
z^{n}_{t}=x^{n}_{t}-j^{n}_{t},
$$
and for any process $z_{t}$ on $\Omega^{n}$ denote
(whenever it makes sense)
\begin{eqnarray}
                                              \label{4.17.3}
M^{n}_{t}(z_{\cdot}):=u(t,z _{t})-u(0,z_{0})
-\int_{0}^{t}u_{t}(s,z _{s})\,ds-
\int_{0}^{t}u_{x^{i}}(s,z _{s})\,dB^{ni}_{s}\nonumber
\\
-(1/2)\int_{0}^{t}u_{x^{i}x^{j}}(s,z _{s})
\,d\< m^{n}\>^{ij}_{s} ,
\end{eqnarray}
$$
\rho^{n}_{s}(z_{\cdot},x)=
 u(s,z _{s}+x)-
u(s,z _{s } )-x^{i}u_{x^{i}}(s,z _{s})
-(1/2)x^{i}x^{j}u_{x^{i}x^{j}}(s,z _{s}),
$$
\begin{equation}
                                                  \label{4.17.4}
R^{n}_{t}(z_{\cdot})=\int_{0}^{t}
\int_{|x|\leq1}\rho^{n}_{s}(z_{\cdot},x)\,\nu^{n}(dsdx).
\end{equation}

Notice that,
 by It\^o's formula
 (see Theorem 1, \S3, Ch.~2 of \cite{LS}) the process
$M^{n}_{t}(z^{n}_{\cdot})-R^{n}_{t}(z^{n}_{\cdot})$
is a local $\cF^{n}_{t}$-martingale.
To be more precise Theorem 1, \S3, Ch.~2 of \cite{LS} says that
$$
M^{n}_{t}(z^{n}_{\cdot})-R^{n}_{t}(z^{n}_{\cdot})
=\sum_{0<s\leq t}
\big[u(s,z^{n}_{s})-u(s,z^{n}_{s-})-u_{x^{i}}(s,z^{n}_{s-})
\Delta z^{ni}_{s}\big]
$$
$$
- \int_{0}^{t}\int_{|x|\leq1}\big[u(s,z^{n}_{s}+x)-
u(s,z^{n}_{s } )-x^{i}u_{x^{i}}(s,z^{n}_{s})\big]\,
\nu^{n}(dsdx)
$$
$$
+\int_{0}^{t}u_{x^{i}}(s,x^{n}_{s-})\,dm^{ni}_{s}.
$$
Here the last term is a local martingale as is any
stochastic integral with respect to
a local martingale and the sum of remaining
terms equals
$$
\int_{0}^{t}\int_{|x|\leq1}\big[u(s,z^{n}_{s-}+x)-
u(s,z^{n}_{s-} )-x^{i}u_{x^{i}}(s,z^{n}_{s-})\big]\,
\bar{\mu}(dsdx)
$$
which is the stochastic integral with respect to the
martingale measure
$\bar{\mu}=\mu-\nu$ and thus also is a
local martingale.

Take the $\cF^{n}_{t}$-stopping time $\gamma^{n}$
introduced in (\ref{4.17.2}). Then
$$M^{n}_{t \wedge\gamma^{n}}(z^{n}_{\cdot})
-R^{n}_{t\wedge\gamma^{n}}(z^{n}_{\cdot}) $$ is again a local
martingale.
It turns out that, for each $T\in[0,\infty)$, the trajectories
of
$
M^{n}_{t \wedge\gamma^{n}}(z^{n}_{\cdot})$, $t\in[0,T]$,
 are bounded
and even uniformly in $n$, Indeed, let $r$
be such that $u(t,x)=0$
for $|x|\geq r$. Notice that $z^{n}_{t}=x^{n}_{t}$
for $0\leq t< \gamma^{n}$. Then we find
$$
 \int_{0}^{t\wedge\gamma^{n}}
u_{x^{i}}(s,z^{n}_{s })\,dB^{ni}_{s}
= \int_{0}^{t\wedge\gamma^{n}}
u_{x^{i}}(s,x^{n}_{s })b^{ni}(s,x^{n}_{\cdot}) \,ds,
$$
where
$$
|u_{x^{i}}(s,x^{n}_{s })b^{ni}(s,x^{n}_{\cdot})|=0
$$
if $|x_{s}|\geq r$ (since $u(t,x)=0$ for $|x|\geq r$) and
$$
|u_{x^{i}}(s,x^{n}_{s })b^{ni}(s,x^{n}_{\cdot})|
\leq  L(r,s) \sup_{s,x}|u_{x}(s,x)|
$$
if  $|x_{s}|\leq r$ (see Assumption \ref{assumption 3.8.2}).
Therefore,
$$
\big|\int_{0}^{t\wedge\gamma^{n}}
u_{x^{i}}(s,z^{n}_{s })\,dB^{ni}_{s}\big|
\leq\sup_{s,x}|u_{x}(s,x)|\int_{0}^{t}L(r,s)\,ds.
$$
Similarly one treats the integrals with respect to
$\< m^{n}\>^{ij}_{s}$. As long as $R^{n}_{t }(z^{n}_{\cdot})$
is concerned we notice that, for $|x|\leq1$
and $0\leq t< \gamma^{n}$, we have
$$
|\rho^{n}_{s}(z^{n}_{\cdot},x)|\leq
 N|x|^{3}I_{|z^{n}_{s}|\leq r+1}
=N|x|^{3}I_{|x^{n}_{s}|\leq r+1},
$$
where the constant $N$ can be expressed
in terms of the third-order
derivatives of $u$ only. Therefore,
$$
|R^{n}_{t\wedge\gamma^{n}}(z^{n}_{\cdot})|
\leq N\theta_{r+1,T}^{n},
$$
where $\theta_{r ,T}^{n}$ is introduced in
 Remark \ref{remark 3.13.1}.
By this remark for any $t$ we have
$E|R^{n}_{t\wedge\gamma^{n}}(z^{n}_{\cdot})|\to 0$. It follows that
$E^{n}|R^{n}_{t\wedge\gamma^{n}}(z^{n}_{\cdot})|<\infty$,
 so that the local martingale $M^{n}_{t\wedge\gamma^{n}}(z^{n}_{\cdot})
-R^{n}_{t\wedge\gamma^{n}}(z^{n}_{\cdot})$ is
in fact a martingale.

Hence,
$$
E^{n}f(x^{n}_{t_{1}},...,x^{n}_{t_{m}})
\big[M^{n}_{t\wedge\gamma^{n}}(z^{n}_{\cdot})
-R^{n}_{t\wedge\gamma^{n}}(z^{n}_{\cdot})-
(M^{n}_{s\wedge\gamma^{n}}(z^{n}_{\cdot})
-R^{n}_{s\wedge\gamma^{n}}(z^{n}_{\cdot})) \big]=0.
$$
Since $E^{n}|R^{n}_{t\wedge\gamma^{n}}
(z^{n}_{\cdot})|\to0$, we also have
$$
\lim_{n\to\infty}E^{n}f(x^{n}_{t_{1}},...,x^{n}_{t_{q}})
\big[M^{n}_{t\wedge\gamma^{n}}(z^{n}_{\cdot}) -
M^{n}_{s\wedge\gamma^{n}}(z^{n}_{\cdot})
  \big]=0 .
$$
Furthermore, due to Remark \ref{remark 4.17.1},
 $P(\gamma^{n}\leq T)\to0$
as $n\to\infty$ for each $T\in[0,\infty)$.
In light of this fact
and by virtue of the uniform boundedness of
$M^{n}_{.\wedge\gamma^{n}}(z^{n}_{\cdot})$, we obtain
\begin{equation}
                                        \label{4.16.2}
\lim_{n\to\infty}E^{n}
\big|M^{n}_{t\wedge\gamma^{n}}(z^{n}_{\cdot})-
M^{n}_{s\wedge\gamma^{n}}(z^{n}_{\cdot})\big|
I_{\gamma_{n}\leq t}=0,
\end{equation}
so that
$$
\lim_{n\to\infty}E^{n}f(x^{n}_{t_{1}},...,x^{n}_{t_{q}})
\big[M^{n}_{t }(z^{n}_{\cdot}) - M^{n}_{s }(z^{n}_{\cdot})
  \big]I_{t<\gamma^{n}}=0.
$$
In addition, obviously, $M^{n}_{t }(z^{n}_{\cdot})=
M^{n}_{t }(x^{n}_{\cdot})$ for
$t<\gamma^{n}$ and in the same way as above
one can prove that  the trajectories
of $M^{n}_{t }(x^{n}_{\cdot})$, $t\in[0,T]$,
are uniformly bounded in $n$ for each $T$.
It follows that
 (\ref{4.16.2}) holds with $t,s,x^{n}_{\cdot}$
in place of
$t\wedge\gamma^{n}$, $s\wedge\gamma^{n}$, $z^{n}_{\cdot}$,
respectively. Thus,
$$
\lim_{n\to\infty}E^{n}f(x^{n}_{t_{1}},...,x^{n}_{t_{q}})
\big[M^{n}_{t }(x^{n}_{\cdot}) - M^{n}_{s }(x^{n}_{\cdot})
  \big] =0
$$
which is rewritten as (\ref{4.21.6}).
The asserted boundedness of the integrand in (\ref{4.21.6})
follows easily from the above argument.
 The lemma is proved.

After we have exploited stochastic bases $(\Omega^{n},\cF^{n},
\cF^{n}_{t},t\geq0,P^{n})$, we will pass to
processes defined on the same probability space.
We are going to rely upon two facts. First
  we know from Theorem 1, \S5, Ch.~6 of \cite{LS} that,
owing to Assumption \ref{assumption 3.8.5},
$\bQ$ is concentrated
on the space of {\em continuous\/} $\bR^{d}$-valued functions
defined on $[0,\infty)$. Second,
 remember that if $y^{n}_{\cdot}\to y_{\cdot}$
in $\cD$ and $y_{\cdot}$ is continuous, then
$|y^{n}_{\cdot}-y_{\cdot}|^{*}_{t}\to0$ for any $t<\infty$,
where
$$
y^{*}_{t}:=\sup_{r\leq t}|y_{r}|.
$$
Owing to these facts and Skorokhod's embedding theorem
(see \S 6, Ch.~1 of \cite{Sk}),
we may assume that all the processes $x^{n}_{\cdot}$,  $n
=1,2,...$,
are given on the same probability space and there is a
continuous process $x_{t}$ such that  (a.s.)
\begin{equation}
                                         \label{4.12.2}
\lim_{n\to\infty}\sup_{t\leq T}|x^{n}_{t}-x_{t}|=0\quad
\forall T\in[0,\infty).
\end{equation}

\begin{lemma}
                                            \label{lemma 4.21.2}
Assume that for any $T$
\begin{equation}
                                            \label{3.13.1}
E\int_{0}^{T}I_{G}(t,x_{t})\,dt=0,
\end{equation}
which is certainly true if $G=\emptyset$.
Then the assertion of Theorem \ref{theorem 3.8.1} holds.
\end{lemma}

Proof. As explained before the lemma
we can write $E$ in place of $E^{n}$ in (\ref{4.21.6}).
Then we  insert $I_{x_{p}\not\in G_{p}}$,
 which is harmless due to (\ref{3.13.1}),
 in the integral in
(\ref{4.21.6}) (notice $x_{p}$ and not $x^{n}_{p}$).
Furthermore, we remember the last  assertion of
Lemma \ref{lemma 4.21.1} and
use Assumption \ref{assumption 3.8.3},  (\ref{4.12.2}),
 and the dominated convergence theorem to conclude that
the limit in (\ref{4.21.6})
equals
\begin{eqnarray}
                                               \label{3.14.3}
Ef(x_{t_{1}},...,x_{t_{q}})
\int_{s}^{t}I_{x_{p}\not\in G_{p}}
\big[
 a^{ ij}(p,x _{p})u_{x^{i}x^{j}}(p,x _{p})\nonumber
\\
+b^{ i}(p,x _{p})u_{x^{i}}(p,x _{p})+u_{p}(p,x _{p})\big]\,dp.
\end{eqnarray}

By using (\ref{3.13.1}) again, we obtain that
$$
 E f(x _{t_{1}},...,x _{t_{q}})
\big[u(t,x _{t})-u(s,x _{s})\big]
$$
$$
= E f(x _{t_{1}},...,x _{t_{q}})
\int_{s}^{t}\big[u_{p}(p,x_{p})+ a^{ ij}(p,x _{p})
u_{x^{i}x^{j}}(p, x _{p})
$$
$$
+b^{ i}(p,x _{p})u_{x^{i}}(p,x _{p})\big]\,dp,
$$
for any bounded continuous $f$ and $t_{i}\leq s\leq t$.
The latter just amounts
to saying that the process
$$
u(t,x_{t})-\int_{0}^{t}\big[u_{s}(s,x_{s})+
a^{ij}(s,x_{s})u_{x^{i}x^{j}}(s,x_{s})
+b^{i}(s,x_{s})u_{x^{i}}(s,x_{s})\big]\,ds
$$
is an $\cF^{x}_{t}$-martingale, where $\cF^{x}_{t}$ is
 the $\sigma$-field generated by
$x_{s}$, $s\leq t$.
It only remains to remember the L\'evy-Doob-Stroock-Varadhan
characterization theorem
 (see, for instance, Sec.~4.5 in \cite{SV} or
  Secs.~2.6 and 2.7 in
\cite{IW}). The lemma is proved.

\begin{remark}{\rm
                                          \label{remark 4.21.3}
In the general case the
above proof and Fatou's theorem show that,
if $f$ {\em is nonnegative\/}, then
$$
 E f(x _{t_{1}},...,x _{t_{q}})
\big[u(t,x _{t})-u(s,x _{s})\big]
$$
$$
\leq E f(x _{t_{1}},...,x _{t_{q}})
\int_{s}^{t}I_{x_{p}\not\in G_{p}}\big[u_{p}
(p,x_{p})+ a^{ ij}(p,x _{p})
u_{x^{i}x^{j}}(p, x _{p})
$$
\begin{equation}
                                             \label{4.21.7}
+b^{ i}(p,x _{p})u_{x^{i}}(p,x _{p})\big]\,dp
+I,
\end{equation}
where
\begin{eqnarray}
                                                \label{4.20.1}
I=Ef(x_{t_{1}},...,x_{t_{q}})
\int_{s}^{t}I_{x_{p}\in G_{p}}\varlimsup_{n\to\infty}
\big[ a^{nij}
(p,x^{n}_{\cdot})u_{x^{i}x^{j}}(p,x^{n}_{p})\nonumber
\\
+b^{ni}(p,x^{n}_{\cdot})u_{x^{i}}(p,x^{n}_{p})
+u_{p}(p,x^{n}_{p})\big]\,dp.
\end{eqnarray}
}
\end{remark}

In the following lemma
 we  finish proving  Theorem \ref{theorem 3.8.1}.
At this moment we take  Theorem \ref{theorem 3.14.1}
for granted.

\begin{lemma}
                                        \label{lemma 4.21.4}
Equation (\ref{3.13.1}) holds and hence,
by Lemma \ref{lemma 4.21.2}, Theorem \ref{theorem 3.8.1}
 holds true as well.
\end{lemma}

Proof. First, we  estimate the   $\varlimsup$ in (\ref{4.20.1}).
Fix $\omega$ and almost any $p$ for which (\ref{4.19.2}) holds
with $p$ in place of $t$ and $x_{p}(\omega)\in G_{p}$. Then we can
replace $\varlimsup_{n\to\infty}$ with
$\lim\limits_{n'\to\infty}$, where $n'$ is an appropriate sequence
tending to infinity. By extracting further subsequences when
necessary we may assume that $a^{n'}(p,x^{n'}_{\cdot})$ and
$b^{n'}(p,x^{n'}_{\cdot})$ converge to some $\bar{a}$ and
$\bar{b}$. Since   $x_{p}\in G_{p}$ and
$|x^{n}_{\cdot}-x_{\cdot}|^{*}_{p}\to0$, (\ref{4.19.2}) implies
that $\det\bar{a}\geq\delta(p,x_{p})$. In addition,
$$
|\bar{b}|+\Trace\,\bar{a}
\leq L(|x_{p}|+1,p)
$$
 due to Assumption \ref{assumption 3.8.2}.
Combined with $\det\bar{a}\geq\delta(p,x_{p})$
this yields
$$
\bar{a}^{ij}\lambda^{i}\lambda^{j}\geq\delta(p,x_{p})
 L^{-(d-1)}(|x_{p}|+1,p) |\lambda|^{2}=:
\bar{\delta}(p,x_{p})|\lambda|^{2}\geq
\tilde{\delta}(p,x_{p})|\lambda|^{2}
$$
for all $\lambda\in\bR^{d}$, where $\tilde{\delta}
=I_{G}\bar{\delta}$.
Now by replacing $\delta$ with $\tilde{\delta}$ and
 both $K(r,t)$ and $L(r,t)$
with $L(r+1,t)$
in Sec.~\ref{section 3.14.1}, we conclude that
$$
\varlimsup_{n\to\infty}
\big[a^{nij}(p,x^{n}_{\cdot})u_{x^{i}x^{j}}(p,x^{n}_{p})
+b^{ni}(p,x^{n}_{\cdot})u_{x^{i}}(p,x^{n}_{p})
+u_{p}(p,x^{n}_{p})\big]
$$
$$
\leq
u_{p}(p,x_{p})+F(p,x _{p},u_{x x }(p,x_{p}))
+ L(|x _{p}|+1,p)|u_{x}(p,x_{p})|.
$$

Furthermore,
Remark \ref{remark 3.13.2} shows that the same estimate
holds for the expression in brackets in (\ref{3.14.3}),
so that according to (\ref{4.21.7})
$$
Ef(x_{t_{1}},...,x_{t_{q}})
\big[u(t,x_{t})-u(s,x_{s})\big]
\leq Ef(x_{t_{1}},...,x_{t_{q}})\int_{s}^{t}\big[
u_{p}(p,x_{p})
$$
$$
+F(p,x _{p},u_{x x }(p,x_{p}))
+ L(|x _{p}|+1,p)|u_{x}(p,x_{p})\big]\,dp
$$
if $f\geq0$.
Hence the process
$$
u(t,x_{t})-\int_{0}^{t}\big[
u_{s}(s,x_{s})
+F(s,x _{s},u_{x x }(s,x_{s}))
+ L(|x _{s}|+1,s)|u_{x}(s,x_{s})\big]\,ds
$$
is a supermartingale and by Theorem \ref{theorem 3.14.1}
estimate (\ref{2.26.2}) holds. If we take there $f=I_{G}$
and remember that the Lebesgue measure of $G$ is zero
and $\bar{\delta}(t,x)>0$ on $G_{t}$
for almost all $t$, then we come to (\ref{3.13.1})
with $T\wedge\tau_{r}$ in place of $T$. Upon letting
$r\to\infty$ we finally obtain (\ref{3.13.1})
as  is. The lemma is proved.

The following theorem is used in Remark \ref{remark 4.25.1}.
Its proof is obtained by changing variables.
We introduce an assumption different from
Assumption~\ref{assumption 3.8.4}.

\begin{assumption}{\rm
                                  \label{assumption 4.25.2}
If $G\ne\emptyset$, then
$G=\bigcup_{m=1}^{\infty}G^{m}$, where $G^{m}$ are Borel
sets. For each $m$,
we are given
  an
integer $d_{m}\geq1$,
 a nonnegative Borel function $\delta_{m} $ defined on $(0,\infty)
\times\bR^{d_{m}}$,
 and a continuous $\bR^{d_{m}}$-valued
function  $v^{m}(t,x)=(v^{m1}(t,x),
 ...,v^{md_{m}}(t,x))$ defined on $[0,\infty)\times\bR^{d}$
and having there continuous in $(t,x)$ derivatives
$v^{mi}_{t},v^{mi}_{x},v^{mi}_{xx}$.
For each $m$ and almost
every $t\in(0,\infty)$,

(i) the set $v^{m}(t,G^{m}_{t})$ has $d_{m}$-dimensional
Lebesgue measure zero,

(ii) for every $x\in v^{m}(t,G^{m}_{t})$
and each sequence  $y^{n}_{\cdot}\in\cD$,
which converges to a continuous function $y_{\cdot}$
 satisfying
  $v^{m}(t,y _{t})=x$,
we have
\begin{equation}
                                               \label{4.25.1}
\varliminf_{n\to\infty}\det V^{mn}(t,y^{n}_{\cdot})
\geq\delta_{m}(t,x)>0,
\end{equation}
where the matrix $V^{mn}(t,y _{\cdot})$ is defined according to
$$
V^{mn}_{ij}(t,y _{\cdot})=v^{mi}_{x^{k}}(t,y_{t})v^{mj}_{x^{r}}
(t,y_{t})a^{nkr}(t,y_{\cdot})\quad i,j=1,...,d_{m}.
$$
}\end{assumption}

\begin{remark}{\rm
Assumption \ref{assumption 3.8.4}
is stronger than Assumption \ref{assumption 4.25.2}.
Indeed, if the former is satisfied, one can take
$G^{m}=G$, $\delta_{m}(t,x)=\delta(t,x)$, $d_{m}=d$,
and $v^{mi}=x^{i}$, $i=1,...,d$, in which case
$\det\,V^{mn}=\det\,a^{n}$.
}\end{remark}

\begin{remark}{\rm

Another case is when
again everything is independent of $m$, but $d_{m}=1$  and
 $v(t,x)=x^{1}$. Then condition (\ref{4.25.1})
becomes
$$
\varliminf_{n\to\infty} a^{n11}(t,y^{n}_{\cdot})
\geq\delta(t,x)>0,
$$
which is much weaker than (\ref{4.19.2}). However, in this case
in order to satisfy requirement (i) of Assumption
\ref{assumption 4.25.2} we need to assume that $G_{t}$
lies in a hyperplane orthogonal to the first coordinate
axis.
}\end{remark}

\begin{remark}{\rm
                                         \label{remark 4.19.4}
Assume that $G=\bigcup_{m=1}^{\infty}G^{m}$, where
 $G^{m}_{t}$ are independent of $t$ and
are hyperplanes $G^{m}_{t}=\{x:(x,\alpha_{m})=\beta_{m}\}$
with certain $\alpha_{m}\in\bR^{d}$
and $\beta_{m}\in\bR$ satisfying
 $|\alpha_{m}|=1$. Assume that we have a Borel
nonnegative functions $\delta_{m}(t,x)$, $x\in\bR$.
Finally, assume that
for every
$m\geq1,t>0$, $x\in\bR^{d}$
such that
$$
(x,\alpha_{m})=\beta_{m},
$$
and each sequence  $y^{n}_{\cdot}\in\cD$,
which converges to a continuous function $y_{\cdot}$ satisfying
  $y _{t}=x$,
we have
$$
\varliminf_{n\to\infty} a^{nij}(t,y^{n}_{\cdot})
\alpha^{i}\alpha^{j}
\geq\delta(t,\beta_{m})>0.
$$

 Then it turns out that Assumption
\ref{assumption 4.25.2} is satisfied. To show this,
it suffices to take
$d_{m}=1$  and $v^{m}(t,x)= (x,\alpha_{m})$ and notice that
the image of $G^{m}_{t}$ under the mapping $v^{m}(t,\cdot)
:G^{m}_{t}\to\bR$ is just one point $\beta_{m}$.
We will use this fact
in Sec.~\ref{section 4.18.1}.

}\end{remark}
\begin{remark}{\rm
Generally, condition (\ref{4.25.1}) is aimed at situations
in which $x^{n}_{t}$ in the limit may
degenerate in some directions
but not along all those which are transversal to $G_{t}$.
}\end{remark}

\begin{theorem}
                                          \label{theorem 4.25.1}

Suppose that Assumptions \ref{assumption 3.8.2},
\ref{assumption 3.8.3},
 \ref{assumption 3.8.5}, and \ref{assumption 4.25.2}
are satisfied and the sequence of distributions
 $(\bQ^{n})_{n\geq1}$ of
$x^{n}_{\cdot}$ converges weakly on $\cD$ to a measure $\bQ$.
Then the assertion of Theorem \ref{theorem 3.8.1}
holds true again.
\end{theorem}

Proof.
We mimic the argument from the proof of
Lemma \ref{lemma 4.21.4}
to show that (\ref{3.13.1}) holds if
  Assumption \ref{assumption 4.25.2} rather than
Assumption \ref{assumption 3.8.4}
is satisfied. The main idea is to change variables
according to the mappings $v^{m}$.

It suffices to prove that, for each $m$,
 equation (\ref{3.13.1})
holds with $G^{m}$ in place of $G$.
Furthermore, without losing generality we may assume that
each set $G^{m}$ is bounded otherwise we could
split each of them into the union of bounded sets and consider
them as new $G^{m}$'s. We fix
 $m,T$, and $R$ and assume that $G^{m}\subset[0,T]\times B_{R}$.
Then the behavior of $v^{m}(t,x)$
for large $|x|$ becomes irrelevant and,
changing $v^{m}$ outside of $[0,T]\times B_{R}$
if necessary, we assume that
\begin{equation}
                                             \label{4.21.2}
v^{m}(t,x) =e_{1}|x|
\end{equation}
for $(t,x)\not\in[0,2T]\times B_{2R}$ , where $e_{1}$
is the first basis vector in $\bR^{d_{m}}$.
It follows that
there is a constant $N_{0}<\infty$ such that
\begin{equation}
                                             \label{4.21.1}
|v^{m}_{x}(t,x)|+|v^{m}_{x^{i}x^{j}}(t,x)|
+|v^{m }_{t}(t,x)|\leq N_{0}\quad\forall t,x.
\end{equation}
It also follows that, for any $r\geq0$,
\begin{equation}
                                               \label{4.21.4}
|v^{m}_{x}(t,x)|\leq r\Longrightarrow |x|\leq 2R+r.
\end{equation}

After that we go back to Lemma \ref{lemma 4.21.1} and take there
$$
u(t,x)=w(t,v^{m}(t,x)),
$$
with $w$ being a function of class
$C^{\infty}_{0}(\bR^{d_{m}+1})$.
By the way,
our stipulation (\ref{4.21.2}) about the
 behavior of $v^{m}$ for large $|x|$
yields that $u\in C^{\infty}_{0}(\bR^{d_{m}+1})$.
We also take the function $f$ in the form
$$
f(y_{1},...,y_{q})=
g(v^{m}(t_{1},y_{1}),...,v^{m}(t_{q},y_{q})),
$$
where $y_{i}\in\bR^{d}$ and $g$ is a continuous
bounded {\em nonnegative\/} function
on $\bR^{qd_{m}}$.  Finally, we define
$$
\tilde{x}^{n}_{t}=v^{m}(t,x^{n}_{t}),\quad
\tilde{x} _{t}=v^{m}(t,x _{t}).
$$
Notice that
$$
a^{nij}(p,x^{n}_{\cdot})u_{x^{i}x^{j}}(p,x^{n}_{p})
+b^{ni}(p,x^{n}_{\cdot})u_{x^{i}}(p,x^{n}_{p})
+u_{p}(p,x^{n}_{p})
$$
$$
=\tilde{a}^{nkr}(p,x^{n}_{\cdot})w_{x^{k}x^{r}}
(p,\tilde{x}^{n}_{p})
+\tilde{b}^{nk}(p,x^{n}_{\cdot})w_{x^{i}}(p,
\tilde{x}^{n}_{p})
+w_{p}(p,\tilde{x}^{n}_{p}),
$$
where, for $y_{\cdot}\in\cD$,
$$
\tilde{a}^{nkr}(p,y_{\cdot})=a^{nij}
(p,y_{\cdot})v^{mk}_{x^{i}}
(p,y_{p})v^{mr}_{x^{j}}(p,y_{p}),
$$
$$
\tilde{b}^{nk}(p,y_{\cdot})=a^{nij}(p,y_{\cdot})
v^{mk}_{x^{i}x^{j}}
(p,y_{p})+b^{ni}(p,y_{\cdot})v^{mk}_{x^{i}}
(p,y_{p})+v^{mk}_{p}(p,y_{p}).
$$
Then on the basis of Fatou's theorem and
 Lemma \ref{lemma 4.21.1}
 we get
$$
E f(x _{t_{1}},...,x _{t_{q}})
\big[u(t,x _{t})-u(s,x _{s})\big]
$$
$$
=E f(x _{t_{1}},...,x _{t_{q}})
\int_{s}^{t}\varlimsup_{n\to\infty}
\big[u_{p}(p,x _{p})+a^{nij}(p,x^{n}_{\cdot})
u_{x^{i}x^{j}}(p,x^{n}_{p})
$$
$$
+b^{ni}(p,x^{n}_{\cdot})u_{x^{i}}(p,x^{n}_{p})\big]\,dp
$$
$$
=Eg(\tilde{x}_{t_{1}},...,\tilde{x}_{t_{q}})
\int_{s}^{t}\varlimsup_{n\to\infty}
\big[w_{p}(p,\tilde{x}_{p})+\tilde{a}^{nkr}(p,x^{n}_{\cdot})
w_{x^{k}x^{r}}(p,\tilde{x}^{n}_{p})
$$
$$
+\tilde{b}^{nk}(p,
x^{n}_{\cdot})w_{x^{k}}(p,\tilde{x}^{n}_{p})\big]\,dp.
$$

Also notice that owing to (\ref{4.21.1}),
$\tilde{a}$ and $\tilde{b}$ satisfy
(\ref{4.14.3}) with $L(r,t)$ replaced with
$N_{0}L(r,t)$. In light of (\ref{4.21.4})  this implies
$$
|\tilde{b}^{n}(t,x^{n}_{\cdot})|
+\Trace\,\tilde{a}^{n}(t,x^{n}_{\cdot})\leq
N_{0}L(2R+|\tilde{x}^{n}_{t}|,t).
$$

In addition, according to
(\ref{4.25.1}), for almost any $t$,
for every $\tilde{x}\in v^{m}(t,G^{m}_{t})$
and each sequence  $y^{n}_{\cdot}\in\cD$,
which converges to a continuous function $y_{\cdot}$
satisfying
  $v^{m}(t,y _{t})=\tilde{x}$,
we have
$$
\varliminf_{n\to\infty}\det \tilde{a}^{n}(t,y^{n}_{\cdot})
\geq\delta_{m}(t,\tilde{x})  >0,
$$
$$
\varliminf_{n\to\infty} \tilde{a}^{nkr}(t,y^{n}_{\cdot})
\lambda^{k}\lambda^{r}
\geq\tilde{\delta}_{m}(t,\tilde{x})|\lambda|^{2}
$$
for all $\lambda\in R^{d_{m}}$, where
$$
\tilde{\delta}_{m}(t,\tilde{x})=
\delta_{m}(t,\tilde{x})
 L^{-(d_{m}-1)}(2R+|\tilde{x}|+1,t)I_{ v^{m}
(G^{m})}(t,\tilde{x})
$$

Then as in the proof of Lemma \ref{lemma 4.21.4}
we find that
$$
Eg(\tilde{x}_{t_{1}},...,\tilde{x}_{t_{q}})
\big[w(t,\tilde{x}_{t})-w(s,\tilde{x}_{s})\big]
\leq Eg(\tilde{x}_{t_{1}},...,\tilde{x}_{t_{q}})
\int_{s}^{t}\big[
w_{p}(p,\tilde{x}_{p})
$$
$$
+F(p,\tilde{x} _{p},w_{x x }(p,\tilde{x}_{p}))
+ L(2R+|\tilde{x}_{p}|+1,p)|w_{x}(p,
\tilde{x}_{p})\big]\,dp,
$$
where the operator $F$ is constructed
on the basis of $\tilde{\delta}_{m}$ and
$N_{0}L(2R+r,t)$
in place of $\delta$ and both $L,K$  from
Sec.~\ref{section 3.14.1},
respectively,
on the space of functions on $\bR^{d_{m}}$
in place of $\bR^{d}$.
Again as in the proof of Lemma \ref{lemma 4.21.4}
 we conclude
that, for any $S $ we have
$$
E\int_{0}^{S}I_{v^{m}(G^{m})}(t,\tilde{x}_{t})\,dt=0.
$$
Since, obviously, $ I_{G^{m}}(t,x)\leq
I_{v^{m}(G^{m})}(t,v^{m}(t,x))$ we get that
(\ref{3.13.1}) holds with $G^{m}$
in place of $G$. As we have pointed out in
the beginning of the proof, this is exactly
what we need. The theorem is proved.

\section{A sufficient condition for precompactness}
                                    \label{section 4.14.2}

One of the conditions of Theorem \ref{theorem 3.8.1}
is that
 the sequence of distributions $(\bQ^{n})_{n\geq1}$ of
$x^{n}_{\cdot}$ on $\cD$ converge.
One can always extract a convergent subsequence
 from a sequence which is   precompact and
here we want to give a simple
sufficient condition for precompactness to hold.
The assumptions of this section are somewhat different
 from the ones of Sec.~\ref{section 4.14.1}
and this was the reason to treat the issue
in a separate section.
We take the objects introduced in
Sec.~\ref{section 4.14.1} before
Assumption \ref{assumption 3.8.2}
and instead of that assumption we require the following.

\begin{assumption}{\rm
                                   \label{assumption 4.14.2}

Assumption \ref{assumption 3.8.2}  is satisfied
with condition (ii) replaced by the following weaker
condition:
 For each $r\in[0,\infty)$ there exists a locally integrable
function $L(r,t)$ given on $[0,\infty)$ such that
$L(r,t)$ increases in $r$ and
$$
|b^{n}(t,y_{\cdot})|+\Trace\,a^{n}(t,y_{\cdot})\leq L(r,t)
$$
whenever $t>0$, $y_{\cdot}\in\cD$, and
$\sup_{s\leq t}|y_{s}|\leq r$.

}\end{assumption}

\begin{lemma}
                                          \label{lemma 4.15.1}
Under Assumptions \ref{assumption 3.8.5}
and \ref{assumption 4.14.2} suppose that
we are given $\cF^{n}_{t}$ stopping times $\tau^{n}_{r}$,
$n=1,2,...,r>0$, and a finite function $\alpha(r)$
defined on $(0,\infty)$
such that we have (i)  for all $n$
and $r$,
\begin{equation}
                                              \label{4.17.6}
|x^{n }_{t}|\leq\alpha(r)\quad\text{if}
\quad 0\leq t<\tau^{n}_{r},
\end{equation}
  and (ii)
\begin{equation}
                                           \label{4.17.7}
\lim_{r\to\infty}\varlimsup_{n\to\infty}P^{n}(\tau^{n}_{r}
\leq T)=0\quad\forall T\in[0,\infty).
\end{equation}
 Then the sequence $(\bQ^{n})_{n\geq1}$
is precompact.

\end{lemma}

Proof. Define
$$
G^{n}_{t}=\int_{0}^{t }
\big[|b^{n}(s,x^{n}_{\cdot})|+\Trace\,a^{n}
(s,x^{n}_{\cdot})\big]\,ds,
$$
$$
F^{n}_{t}=G^{n}_{t}+ \int_{0}^{t }
\int_{|x|>1}\nu^{n}(dsdx)
$$
Owing to Assumption  \ref{assumption 3.8.5},
by Theorem VI.4.18 and Remark VI.4.20 of \cite{JS}
to prove the theorem it suffices to check that
the sequence of distributions on $\cD$
of   $F^{n}_{\cdot}$ is $C$-tight, that is
precompact and each limit point of this
sequence is the distribution of a continuous
process. In turn, due to Theorem VI.4.5
and Remark VI.4.6 (3)  of \cite{JS}, to prove the $C$-tightness
it suffices to show that, for any $T\in[0,\infty)$
and $\varepsilon>0$,
 \begin{eqnarray}
                                        \label{4.17.1}
  &\lim_{N\to\infty}\varlimsup_n&P^n\Big(\sup_{t\le
T}\big|F^{n}_t\big|\ge N\Big)=0,\nonumber
\\
 &\lim_{\delta\downarrow 0}\,\varlimsup_n \, &
P^n\Big(\sup_{t+s\leq T,0\leq s\leq\delta}
\big|F^{n}_{t+s}-F^{n}_t\big|\ge
\varepsilon\Big)=0.
\end{eqnarray}

In view of
Assumption  \ref{assumption 3.8.5} we need only prove
(\ref{4.17.1}) for $G^{n}$ in place of $F^{n}$. We do this
replacement and after that notice that, for any $r$, the left-hand
side of the first equation in (\ref{4.17.1}) is less than
$$
\lim_{N\to\infty}\varlimsup_nP^n\Big(\sup_{t\le
T\wedge\tau^{n}_{r}}\big|G^{n}_t\big|\ge N\Big)
+\varlimsup_{n\to\infty}P^{n}(\tau^{n}_{r} \leq T).
$$
Here the first term is zero for each $r$ since
$G^{n}_t$ is continuous in $t$ and
$$
|G^{n}_t|\leq\int_{0}^{t}L(\alpha(r),u)\,du
$$
for $t<\tau^{n}_{r}$ when by our assumptions $|x^{n}_{t }|\leq r$.
In addition, the second term can be made as small as we wish by
choosing a sufficiently large $r$. This proves the first equation
in (\ref{4.17.1}).

Similarly, the left-hand side of the second equation
in (\ref{4.17.1}) with $G^{n}$ in place of $F^{n}$
is less than
$$
\lim_{\delta\downarrow 0}\varlimsup_nP^n\Big(\sup_{t+s\le
T\wedge\tau^{n}_{r},0\le
s\le\delta}\big|G^{n}_{t+s}-G^{n}_{t}\big|\ge\varepsilon\Big)
+\varlimsup_{n\to\infty}P^{n}(\tau^{n}_{r} \leq T),
$$
where again the first term vanishes since
$$
|G^{n}_{t+s}-G^{n}_{t}|\leq\int_{t}^{t+s}L(\alpha(r),u)\,du.
$$
 The
lemma is proved.

\begin{remark}{\rm
It may be worth noticing that the combination of assumptions (i) and (ii)
of Lemma \ref{lemma 4.15.1} is equivalent to the following:
for any $T\in(0,\infty)$, the sequence of distributions
of   $\sup_{t\le T}|x^n_t|$
is tight or put otherwise
$$
\lim_{r\to\infty}\varlimsup_{n\to\infty}P^n(\sup_{t\le T}|x^n_t|\geq
r)=0.
$$
}\end{remark}
Lemma \ref{lemma 4.15.1} reduces the investigation of
precompactness to estimating $|x^{n}|^{*}_{t}$. Here the following
coercivity assumption turns out to be useful.
\begin{assumption}{\rm
                                  \label{assumption 4.14.1}
For any $n$, there exists a
 nonnegative $\cF^{n}_{t}$-predictable function $L_{n}(t)$ such that

\begin{equation}
                                           \label{4.14.5}
  b^{ni}(t,x^{n}_{\cdot})x^{ni}_{t}
+\Trace\, a^{n}(t,x^{n}_{\cdot})\leq L_{n}(t)(1+|x^{n}_{t}|^{2})
\end{equation}
for almost all $(\omega,t)$. Furthermore,
for any $T\in[0,\infty)$,
$$
\lim_{c\to\infty}
\varlimsup_{n\to\infty}P^{n}\big(\int_{0}^{T}L_{n}(t)\,dt
>c)=0.
$$
}\end{assumption}

\begin{remark}{\rm
                                      \label{remark 4.25.4}
Quite often one imposes a linear growth
assumptions on the  coefficients $a^{n}$ and $b^{n}$,
which of course implies (\ref{4.14.5}).
However, say
in one dimension, if $a^{n}\equiv0$ and
$b^{ni}(t,y_{\cdot})=b^{n}(t,y_{t})$
and $b^{n}(t,y_{t})\geq0$
for $y_{t}<0$ and $b^{n}(t,y_{t})\leq0$
for $y_{t}>0$, then (\ref{4.14.5}) is satisfied
with $L\equiv0$. Therefore generally (\ref{4.14.5})
does not provide any control on the behavior
of $|b^{n}(t,y_{t})|$ for large $|y_{t}|$.

For that reason, Theorem \ref{theorem 4.17.1} below
does not follow from the results of \cite{JS}
and \cite{LS}.
}\end{remark}

\begin{theorem}
                                   \label{theorem 4.17.1}
Let
\begin{equation}
                                       \label{4.17.9}
\lim_{N\to\infty}\varlimsup_{n\to\infty}
P^{n}(|x^{n}_{0}|\geq N)=0
\end{equation}
and let Assumptions \ref{assumption 3.8.5},
 \ref{assumption 4.14.2}, and
\ref{assumption 4.14.1} be satisfied. Then
the sequence $(\bQ^{n})_{n\geq1}$ is precompact.
Furthermore, let $k$ be an integer and $f^{n}(t,x)$ be
  Borel
$\bR^{k}$-valued functions   defined on
$(0,\infty)\times\bR^{d}$
such that $|f^{n}(t,x)|\leq L(|x|,t)$ for all $t,x,n$.
Define
$$
y^{n}_{t}=\int_{0}^{t}f^{n}(s,x^{n}_{s})\,ds.
$$
Then the sequence of distributions of
$(x^{n}_{\cdot},y^{n}_{\cdot})$
on $\cD([0,\infty),\bR^{d+k})$ is precompact as well.

\end{theorem}

Proof. We are going to use a method introduced
in Sec.~4, Ch.~II of \cite{KR}.
Define
$$
z^{n}_{t}=x^{n}_{t}-j^{n}_{t},\
\phi_n(t)=\exp \Big( -2\int_{0}^{t}L
_n(s)\,ds \Big) ,
\ u_n(t,x)=(1+|x|^{2})\phi_n(t).
$$
Also
as in the proof of Lemma \ref{lemma 4.21.1},
use notation (\ref{4.17.3}) and (\ref{4.17.4})
and notice that due to special choice of $u$,
we have $R^{n}_{t}(z_{\cdot})\equiv0$.

Then by  using It\^o's
formula, we get that the process
$$
M^{n}_{t}:=(1+|z^{n}_{t}|^{2})\phi_n(t)-
(1+|x^{n}_{0}|^{2})
$$
\begin{equation}
                                         \label{4.17.5}
-
\int_{0}^{t}\big[2z^{ni}_{s}b^{ni}(s,x^{n}_{\cdot})
+2\Trace\,a^{n}(s,x^{n}_{\cdot})
-2L_n(s)(1+|z^{n}_{s}|^{2})\big]
\phi_n(s)\,ds
\end{equation}
is a local martingale.

Now
take $\gamma^{n}$ again from (\ref{4.17.2})
and remember that $z^{n}_{s}=x^{n}_{s}$ for $s<\gamma^{n}$,
so that the expression in the brackets in (\ref{4.17.5})
is negative due to Assumption \ref{assumption 4.14.1}.
Then we see that
$$
H^{n}_{t}:=(1+|z^{n}_{t\wedge\gamma^{n}}
|^{2})\phi_n(t\wedge\gamma^{n})-
(1+|x^{n}_{0}|^{2})
$$
is a local supermartingale. For any constant $N>0$, the process
$H^{n}_{t}I_{|x^{n}_{0}|\leq N}$ also is
 a local supermartingale
and, since it is bounded from below by the constant
 $-(1+N^{2})$, it is a supermartingale.
Therefore, upon defining
$$
\kappa^{n}_{r}=\inf\{t\geq0:\sup_{s\leq t}|x^{n}_{s}|>r \},
\quad\tau^{n}_{r}=\gamma^{n}\wedge\kappa^{n}_{r},
$$
 we get that, for any $T\in[0,\infty)$,
$$
E^{n}\big(1+|z^{n}_{T\wedge\tau^{n}_{r}}|^{2}\big)
\phi_n(T\wedge\tau^{n}_{r})I_{|x^{n}_{0}|\leq N}
\leq1+N^{2},
$$
$$
E^{n}\big(1+|z^{n}_{ \tau^{n}_{r}}|^{2}\big)
\phi_n( \tau^{n}_{r})I_{|x^{n}_{0}|\leq N,
\tau^{n}_{r}\leq T<\gamma^{n}}
\leq1+N^{2}.
$$
Then we notice that on the interval $[0,\gamma^{n})$
the process $j^{n}_{t}$ is identically zero. Hence,
for $\tau^{n}_{r}\leq T<\gamma^{n}$ we have
$$
|z^{n}_{ \tau^{n}_{r}}|=|x^{n}_{\tau^{n}_{r}}|
=|x^{n}_{\kappa^{n}_{r}}|\geq r
$$
and we obtain
$$
e^{-c}(1+r^{2})P^{n}\bigg(\int_{0}^{T}L_{n}(t)\,dt
\leq c,|x^{n}_{0}|\leq N,
\tau^{n}_{r}\leq T<\gamma^{n}\bigg)\leq1+N^{2},
$$
$$
\lim_{r\to\infty}\varlimsup_{n\to\infty}
P^{n}(|x^{n}_{0}|\leq N,
\tau^{n}_{r}\leq T<\gamma^{n})=0.
$$

This holds for any $N$ and along with
assumption  (\ref{4.17.9})
and Remark \ref{remark 4.17.1} leads first to
to
$$
\lim_{r\to\infty}\varlimsup_{n\to\infty}
P^{n}(
\tau^{n}_{r}\leq T<\gamma^{n})=0
$$
and then to (\ref{4.17.7}).

Finally,
observe that (\ref{4.17.6}) is obviously satisfied
even if $0\leq t<\kappa^{n}_{r}$ rather than
$0\leq t<\tau^{n}_{r}$. Hence, by referring to
Lemma \ref{lemma 4.15.1} we finish proving
the assertion of our theorem regarding the distributions
of $x^{n}_{\cdot}$.

Lemma \ref{lemma 4.15.1} yields the result for
$(x^{n}_{\cdot},y^{n}_{\cdot})$ as well since, obviously,
  for $0\leq t<r\wedge\tau^{n}_{r}$, we have
$$
|y^{n}_{t}|\leq\int_{0}^{r}L(r,s)\,ds.
$$
The theorem is proved.

\section{An example of queueing model}
                                    \label{section 4.18.1}

We consider
a particular queueing system with $d$ service stations
and $d+1$ incoming streams of customers.
We refer the reader to \cite{FS}
for relations of this system to practical problems.
The first $d$ streams are composed of
  customers ``having appointments'', meaning that the customers
from the $i$th stream only go to the $i$th
service station. The last stream, to which
we assign number 0, is the one of ``free''
customers who, upon ``checking in'',
are routed to the service stations according to certain rule
to be described later.  We assume that each service station
consists of infinitely many servers, so that infinitely
many customers
can be served at each station simultaneously.
Denote by $Q^{i}_{t}$
the number of customers
being served at the $i$th station at time $t$.

With  station $i$, $i=1,...,d$, we associate a ``cost''
$\alpha_{i}>0$ and suppose that
  a ``free'' customer arriving at time $t$
is directed to the $i$th station if $i$ is the smallest
integer satisfying
$$
\alpha_{i}Q^{i}_{t-}\leq\alpha_{j}Q^{j}_{t-}\quad
\text{for all}\quad j\ne i.
$$
Such a routing policy is called load-balancing in
\cite{FS}. Here and below in this section
the summation convention over repeated indices {\em
is not enforced\/}.

We take some numbers $\lambda_{0},...,\lambda_{d}>0$
and assume that the $i$th stream
of customers forms a Poisson process with parameter
$\lambda_{i}$.  To describe the service times
we fix some ``thresholds'' $N^{1},...,N^{d}$,
which are positive integers,  and assume that,
given $0<Q^{i}_{t}< N^{i}$, each of $ Q^{i}_{t}$
customers at the $i$th station

(i) has its own server,

(ii) spends
with its server a random time having
exponential distribution with parameter 1,

(iii)  after having been served leaves the system.

However,
given $Q^{i}_{t}\geq N^{i}$, the service
is organized differently. All $Q^{i}_{t}$ customers
are divided into disjoint groups
each consisting of two persons apart from at most
one group having only one member.
Then each of those groups is supposed
to get service according to the rules (i)-(iii)
above. By the way, it is not hard to understand
that on average both discipline of servicing
yield the same
number of customers having been served during
one unit of time.

Finally, we assume that all service times
and arrival processes are as independent
as they can   be.

Now we describe the model in rigorous terms.
For any numbers $y^{1},...,y^{d}$ define
$$
\argmin_{k=1,...,d}y^{k}=i
$$
if $i$ is the least of $1,...,d$ such that $y^{i}\leq y^{k}$
for $k\ne i$. For $x\in\bR^{d}$ and $i=1,...,d$, let
$$
\delta^{i}(x)=\left\{\begin{array}{ll}1&\quad\text{if}\quad
i=\argmin\limits_{k=1,...,d}\alpha_{k}x^{k},
\\
0 &\quad\text{otherwise}.
\end{array}\right.
$$
Take independent Poisson processes $\Pi^0_t,...,
\Pi^d_t$ with parameters $\lambda_{0},...,\lambda_{d}$,
respectively. Then we think of the number of arrivals
at the $i$th station as given by
$$
A^{i}_{t}=\int_{0}^{t}\delta^{i}(Q_{s-})\,d\Pi^{0}_{s}
+\Pi^{i}_{t},
$$
where $Q_{s}=(Q^{1}_{s},...,Q^{d}_{s})$
and $Q^{i}_{t}$ are some integer-valued right continuous
processes having left limits.
To model the number of
departures $D^{i}_{t}$ from the $i$th station up to time $t$
we take  Poisson processes $\Pi^{ij} _{t}$ and
$\Lambda^{ij} _{t}$, $i=1,...,d$, $j=1,2,...$,
having parameter 1 and mutually independent
and independent of
$(\Pi^{0}_{\cdot},...,\Pi^{d}_{\cdot})$.
Then we define
$$
D^{i}(t)=\int_{0}^{t}I_{N^{i}> Q^{i}_{s-}}\sum_{j\geq1}
I_{ Q^{i}_{s-}\geq j}\,d\Pi^{ij}_{ s}
$$
$$
+  \int_{0}^{t}I_{N^{i}\leq Q^{i}_{s-}}\sum_{j\geq1}
\big(I_{
Q^{i}_{s-}\geq2j}+I_{
Q^{i}_{s-}+1\geq2j}\big)\, d\Lambda^{ij}_{ s}.
$$
To be consistent with the description,
$Q_{t}$ should satisfy
the balance equations $Q^{i}_{t}=Q^{i}_{0}+
A^{i}_{t}-D^{i}_{t}$.
Thus, we are going to investigate the system
of equations
$$
dQ^{i}_{t}=\delta^{i}(Q_{t-})\,d\Pi^{0}_{t}
+d\Pi^{i}_{t}-
I_{N^{i}> Q^{i}_{t-}}\sum_{j\geq1}
I_{ Q^{i}_{t-}\geq j}\,d\Pi^{ij}_{t}
$$
\begin{equation}
                                              \label{4.18.1}
- I_{N^{i}\leq Q^{i}_{t-}}\sum_{j\geq1}
\big(I_{
Q^{i}_{t-}\geq2j}+I_{
Q^{i}_{t-}+1\geq2j}\big)\, d\Lambda^{ij}_{t}\quad i=1,...,d.
\end{equation}
Needless to say that we assume that all
the Poisson processes we are dealing with are given
on a probability basis satisfying the ``usual'' assumptions.
We also assume that the initial condition $Q_{0}$
is independent of the Poisson processes.

Notice that
for any   initial condition $Q_{0}$
 there is a unique solution
of (\ref{4.18.1}). Indeed obviously, for any solution
we have
$Q^{i}_{t}\leq Q_{0}^{i}+ \Pi^{0}_{t}+\Pi^{i}_{t}$,
so that, while solving (\ref{4.18.1})
for $t\in[0,T]$, one can safely replace the
infinite sums in (\ref{4.18.1}) with
the sums over $j\leq Q_{0}^{i}+ \Pi^{0}_{T}+\Pi^{i}_{T}$.
After that one  solves  (\ref{4.18.1}) on each
$\omega$ noticing that  between the  jumps
of the Poisson processes $Q_{t}$ is constant and
the jumps of $Q_{t}$
themselves are given by (\ref{4.18.1}).

For obvious reasons we rewrite (\ref{4.18.1})
in terms of representation (\ref{4.24.2}).
First, for $k=0,...,d,i=1,...,d,j\geq1$,  we define
$$
\bar{\Pi}^{k}_{t}=\Pi^{k}_{t}-\lambda_{k}t,\quad
\bar{\Pi}^{ij}_{t}=\Pi^{ij}_{t}-t,\quad
\bar{\Lambda}^{ij}_{t}=\Lambda^{ij}_{t}- t.
$$
These processes are square integrable martingales with
$$
\<\bar{\Pi}^{k}\>_{t}=\lambda_{k}t,\quad
\<\bar{\Pi}^{ij}\>_{t}= t,\quad
\<\bar{\Lambda}^{ij}\>_{t}= t.
$$
Next, for $i=1,...,d$, define
$$
M^{ i}_{t}=\int_{0}^{t}
\delta^{i}(Q _{s-})\,d\bar{\Pi}^{0}_{s}
+ \bar{\Pi}^{i}_{t}-\int_{0}^{t}
I_{N^{i}> Q^{i}_{s-}}\sum_{j\geq1}
I_{ Q^{i}_{s-}\geq j}\,d\bar{\Pi}^{ij}_{s}
$$
$$
- \int_{0}^{t}I_{N^{i}\leq Q^{i}_{s-}}\sum_{j\geq1}
\big(I_{
Q^{i}_{s-}\geq2j}+I_{
Q^{i}_{s-}+1\geq2j}\big)\, d\bar{\Lambda}^{ij}_{s},
$$
which are at least   locally square integrable martingales.
Then after observing  that, for any integer  $q\geq0$,
$$
\sum_{j\geq1}I_{q\geq j}=q,\quad
\sum_{j\geq1}(I_{q\geq2j}+I_{q+1\geq2j})=q,
$$
we turn equation (\ref{4.18.1}) into the equation
\begin{equation}
                                          \label{4.24.3}
dQ^{i}_{t}=(\lambda_{0}\delta^{i}(Q_{t})
+\lambda_{i}-Q^{i}_{t})\,dt+dM^{i}_{t}.
\end{equation}

In order to explain what follows
(in no way is this explanation used in the
proof of Theorem \ref{theorem 4.18.1} below), notice
that (\ref{4.24.3}) seems to imply that
\begin{equation}
                                             \label{4.24.1}
(EQ^{i}_{t} )'=\lambda_{0}
E\delta^{i}(Q_{t})+\lambda_{i}-EQ^{i}_{t}.
\end{equation}
We are interested in the behavior of $Q_{t}$ when
$\lambda_{i}$'s are large but $\lambda_{0}$
is much smaller than $\lambda_{1},...,\lambda_{d}$.
 Then, on the one hand, $EQ^{i}_{t}$
should be large for moderate $t$ and, on the other hand,
 the first term
on the right in (\ref{4.24.1}) can be neglected.
In that situation  equation (\ref{4.24.1}) turns out to have
a stable point $EQ_{t}^{i}\equiv\lambda_{i}$.
This means that, if for the initial condition
we have $EQ_{0}^{i}=\lambda_{i}$, then
$EQ_{t}^{i}=\lambda_{i}$ for all $t$.
Notice that since $\lambda_{i}$'s are large,
so should be $EQ_{0}^{i}$.

Therefore, we
write $\lambda_{i}=\bar{\lambda}_{i}+\Delta\lambda_{i}$,
where $\Delta\lambda_{i}$ will be assumed to have order of
$\lambda_{0}$,
denote
$$
\bar{Q}^{i}_{t}=Q^{i}_{t}-\bar{\lambda}_{i}
$$
and rewrite (\ref{4.24.3}) in terms of $\bar{Q}_{t}$.
At this moment we introduce the assumption that
\begin{equation}
                                            \label{4.24.4}
\bar{\lambda}_{i}\alpha_{i}=n,\quad i=1,...,d,
\end{equation}
with $n$ being an integer
(independent of $i$) to be sent to infinity.
This is convenient
due to the simple fact that then
$$
\delta^{i}(x)=\delta^{i}(x-\bar{\lambda}).
$$
In this notation (\ref{4.24.3}) becomes
$$
d\bar{Q}^{i}_{t}= (\lambda_{0}\delta^{i}(\bar{Q}_{t})
+\Delta\lambda_{i} -\bar{Q}^{i}_{t})\,dt+dM^{i}_{t}.
$$
To understand what kind of normalization is natural
we compute the quadratic characteristics
of $M^{i}_{t}$. Notice that,
for any integer $q\geq0$, we have
$$
\sum_{j\geq1}(I_{q\geq2j}+I_{q+1\geq2j})^{2}=
\sum_{j\geq1}(I_{q\geq2j}+2I_{q\geq2j}+I_{q+1\geq2j})
$$
$$
=3[q/2]+[(q+1)/2]=:qf(q),
$$
where $[a]$ is the integer part of $a$. By the way,
we can only define $f(q)$ by
the above formula for all real $q>0$. If $q\leq0$,
 we let $f(q)=0$.
Then
\begin{equation}
                                            \label{4.24.7}
0\leq f\leq2,\quad\lim_{q\to\infty}f(q)=2.
\end{equation}
It follows that
$$
d\<M\>^{ii}_{t}=[\lambda_{0}\delta^{i}(\bar{Q}_{t})
+\lambda_{i}+Q^{i}_{t}I_{Q^{i}_{t}<N^{i}}
+Q^{i}_{t}f(Q^{i}_{t})I_{Q^{i}_{t}\geq N^{i}}]\,dt.
$$
Also due to independence of our Poisson processes
and the fact that $\delta^{i}\delta^{j}=0$
for $i\ne j$, we get
$$
\<M\>^{ij}_{t}=0\quad\text{for}\quad i\ne j.
$$
If we believe that, in a sense, $Q^{i}_{t}\sim\lambda_{i}$,
then $Q^{i}_{t}/\lambda_{i}$ should converge
as well as $M^{i}_{t}/\sqrt{\lambda}_{i}$, and we see
that it is natural to expect $\bar{Q}^{i}_{t}
/\sqrt{\lambda}_{i}$
to converge to certain limit. To make the model
 more meaningful
we also assume that the thresholds $N^{i}$'s are large
and roughly speaking proportional to $\lambda_{i}$.
In this way we convince ourselves that the following result
seems natural.

\begin{theorem}
                                       \label{theorem 4.18.1}
Let $\alpha_{1},...,\alpha_{d}>0$ and
  $\mu_{0},...,\mu_{d}\geq0$ and
$\nu_{1},...,\nu_{d}\in\bR$ be fixed parameters.
For   $n=1,2,...$ define
$$
\lambda_{i}=n\alpha_{i}^{-1}+\mu_{i}\sqrt{n},
\quad i=1,...,d,\quad\lambda_{0}=\mu_{0}\sqrt{n},
$$
$$
N^{i}=n\alpha_{i}^{-1}+\nu_{i}\sqrt{n},\quad i=1,...,d.
$$
Let $Q_{t}=Q_{t}^{n}$ be the solution of
(\ref{4.18.1}) with certain initial condition
independent of the Poisson processes
and introduce
$$
x_{t}^{n}=
n^{-1/2}(Q^{n1}_{t}-n\alpha_{1}^{-1},...
,Q^{nd}_{t}-n\alpha_{d}^{-1}).
$$
Let $\bQ^{n}$ be the distribution
of $x_{\cdot}^{n}$ on $\cD$. Finally, assume that
the distribution of
$
x^{n}_{0}
 $
weakly converges to a distribution $F_{0}$
as $n\to\infty$.

Then, as $n\to\infty$,
$\bQ^{n}$ converges weakly to the distribution
of a solution of the following system
\begin{equation}
                                         \label{4.24.6}
dx^{i}_{t}=(\mu_{0}\delta^{i}(x_{t})+
\mu_{i}-x^{i}_{t})\,dt+
\alpha_{i}^{-1/2}(2
+I_{x^{i}_{t}\geq\nu_{i}})^{1/2}\,dw^{i}_{t}, \quad i=1,...,d
\end{equation}
considered on some probability space with $w_{t}$
being a $d$-dimensional Wiener process and $x_{0}$
distributed according to $F_{0}$.
\end{theorem}

Proof. First of all notice that (\ref{4.24.6})
has   solutions on appropriate probability spaces
 and any  solution
has the same distribution on the space
of $\bR^{d}$-valued continuous functions.
This follows from the fact that an obvious
change of probability measure allows us
to consider the case with no drift terms in
(\ref{4.24.6}). In that case (\ref{4.24.6})
becomes just a collection of unrelated
one-dimensional equations
with uniformly nondegenerate and bounded
diffusion. Weak unique solvability
of such equations is a very well known fact
(see, for instance, Theorems 2 and 3 of \cite{Kr69}).

In the proof of convergence
we will be using Theorems \ref{theorem 4.17.1}
and \ref{theorem 3.8.1}.
Observe that Assumption \ref{assumption 3.8.5}
is satisfied since $x^{ni}_{t}$ has no jumps
bigger than $2n^{-1/2}$ and $\nu^{n}((0,\infty)
\times B^{c}_{a})=0$ if $n>4d/a^{2}$.
Furthermore, if in  the argument before the theorem
 we take
 $\bar{\lambda}_{i}
=n\alpha_{i}^{-1}$, so that
  (\ref{4.24.4}) holds,
and let
$
\Delta\lambda_{i}=\mu_{i}\sqrt{n},
$
then after noticing that, by definition,
$$
Q^{ni}=n^{1/2}x^{ni}_{t}+n\alpha_{i}^{-1},
$$
we easily obtain
\begin{equation}
                                             \label{4.25.2}
dx^{n}_{t}=b^{n}(x^{n}_{t})\,dt+dm^{n}_{t},
\quad\<m^{n}\>_{t}
=\int_{0}^{t}a^{n}(x^{n}_{s})\,ds,
\end{equation}
where
$$
b^{ni}(x)=\mu_{0}\delta^{i}(x)+\mu_{i}-x^{i},\quad
a^{nij}(x)=\delta^{ij}\big(n^{-1/2}
\mu_{0}\delta^{i}(x)+ \alpha_{i}^{-1}
+\mu_{i}n^{-1/2}
$$
$$
+(x^{i}n^{-1/2}+\alpha_{i}^{-1})_{+}\big[I_{x^{i}<\nu^{i}}
+f(n^{1/2}\,x^{i}+n\alpha_{i}^{-1})I_{x^{i}\geq\nu^{i}}\big]
\big).
$$
Upon remembering (\ref{4.24.7}) we see that,
for a constant $N$ and all $n$ and $x$, we have
$|b^{n}(x)|+\Trace\,a^{n}(x)\leq N(1+|x|)$,
which shows that
Assumptions \ref{assumption 3.8.2}
and \ref{assumption 4.14.2},
equivalent in our present situation,
and Assumption \ref{assumption 4.14.1} are satisfied.
By Theorem \ref{theorem 4.17.1} the sequence
$(\bQ^{n})$ is precompact.

Next, obviously Assumption \ref{assumption 3.8.3}
is satisfied if we take
$$
G=\{(t,x):t>0,\prod_{i,j=1}^{d}(\alpha_{i}x^{i}
-\alpha_{j}x^{j})(x^{i}-\nu_{i})=0\},
$$
$$
b^{i}(x)=\mu_{0}\delta^{i}(x)+\mu_{i}-x^{i},\quad
a^{ij}(x)=\delta^{ij} \alpha_{i}^{-1}\big(
1+ I_{x^{i}<\nu^{i}}
+2I_{x^{i}\geq\nu^{i}}\big) .
$$
Finally, Assumption \ref{assumption 3.8.4} is satisfied
since $\det\,a^{n}(x)\geq\alpha_{1}^{-1}
\cdot...\cdot\alpha_{d}^{-1}$
everywhere.

By Theorem \ref{theorem 3.8.1} every convergent subsequence
of $(\bQ^{n})$ converges to the distribution
of a solution of (\ref{4.24.6}) with the above specified
initial distribution. Since all such
solutions have the same distribution, the whole sequence
$(\bQ^{n})$ converges to the distribution
of any solution of (\ref{4.24.6}). The theorem is proved.

\begin{remark}{\rm

                                                \label{remark 5.19.2}
In Theorem \ref{theorem 4.18.1} we assume that
$ Q^{ni}_{0}$
goes to infinity with certain rate, namely
$ Q^{ni}_{0}\sim n\alpha^{-1}_{i} $.
Interestingly enough, if we change the rate,
the diffusion approximation changes. Indeed, keep
all the assumption of Theorem \ref{theorem 4.18.1}
apart from the assumption that $x^{n}_{0}$
converges in distribution and instead assume that,
for a $\gamma\in[0,\infty)$ say for $\gamma=0$,
$$
n^{-1/2}(Q^{n1}_{0}-n\gamma\alpha^{-1}_{1},...,Q^{nd}_{0}-n
\gamma\alpha^{-1}_{d})
$$
 converges in law to a random vector.
Notice that the case $\gamma=1$ is covered by
Theorem \ref{theorem 4.18.1}.
We claim that, for $\gamma>1$, the processes
$$
y^{n}_{t}=n^{-1/2}(Q^{n1}_{t}-nq_{t}\alpha^{-1}_{1},
...,Q^{nd}_{t}-nq_{t}\alpha^{-1}_{d}),
$$
where $q_{t}=1+(\gamma-1)e^{-t}$,
weakly converge to a solution of the system
$$
dy^{i}_{t}=(\mu_{0}\delta^{i}(y_{t})+
\mu_{i}-y^{i}_{t})\,dt+
\alpha_{i}^{-1/2}(1+q_{t})^{1/2}\,dw^{i}_{t},
\quad  i=1,...,d,
$$
and for $\gamma\in[0,1)$ weakly converge to a solution of
$$
dy^{i}_{t}=(\mu_{0}\delta^{i}(y_{t})+
\mu_{i}-y^{i}_{t})\,dt+
\alpha_{i}^{-1/2}(1+2q_{t})^{1/2}\,dw^{i}_{t},
\quad  i=1,...,d.
$$

Indeed,
we have
$$
Q^{ni}=n^{1/2}y^{ni}_{t}+nq_{t}\alpha_{i}^{-1},
\quad dq_{t}=(1-q_{t})\,dt,
$$
$$
dy^{n}_{t}=b^{n}(y^{n}_{t})\,dt+dm^{n}_{t},
\quad\<m^{n}\>_{t}
=\int_{0}^{t}a^{n}(y^{n}_{s})\,ds,
$$
where
$$
b^{ni}(x)=\mu_{0}\delta^{i}(x)+\mu_{i}-x^{i},\quad
a^{nij}(x)=\delta^{ij}\bigg(n^{-1/2}
\mu_{0}\delta^{i}(x)+ \alpha_{i}^{-1}
$$
$$
+\mu_{i}n^{-1/2}
+(x^{i}n^{-1/2}+q_{t}\alpha_{i}^{-1})_{+}\big[I_{
(\gamma-1)e^{-t}<\alpha_{i}(\nu^{i}-
x^{i})n^{-1/2}}
$$
$$
+f(n^{1/2}\,x^{i}+nq_{t}\alpha_{i}^{-1})
I_{
(\gamma-1)e^{-t}\geq\alpha_{i}(\nu^{i}-
x^{i})n^{-1/2}}\big]
\bigg).
$$

As in the proof of Theorem \ref{theorem 4.18.1}
one checks that the sequence of distributions of
$y^{n}_{\cdot}$ is precompact. Furthermore, obviously,
for any $x$
$$
a^{nij}(x)\to\left\{\begin{array}{ll}
\delta^{ij}\alpha_{i}^{-1}(1+q_{t} )&\quad
\text{if}\quad\gamma<1,
\\
\delta^{ij}\alpha_{i}^{-1}(1+2q_{t} )&\quad
\text{if}\quad\gamma>1,
\end{array}\right.
$$
and this yields our claim in the same way as in
the proof of Theorem~\ref{theorem 4.18.1}.

}\end{remark}
\begin{remark}{\rm
                                                \label{remark 5.19.3}
We tried to explain before the proof of
Theorem \ref{theorem 4.18.1} why its statement looks natural.
Now we can also explain how the function $q_{t}$
from Remark \ref{remark 5.19.2} was found.
The explanations is based on a kind of law of large numbers
which in queueing theory is associated with
 so-called ``fluid approximations''.
Generally, ``fluid approximations'' can also be
derived from Theorems \ref{theorem 4.17.1}
and \ref{theorem 3.8.1}. For instance, if
$\lambda_{k}=\lambda_{k}(n)$ and $\lambda_{k}(n)/n
\to\beta_{k}$ as $n\to\infty$, and $\beta_{0}=0$,
then under the condition that
$Q^{n}_{0}/n$ converges in probability to a constant
vector, the processes $Q^{n}_{t}/n$
converge in probability uniformly on each finite time interval
to the deterministic solution
of the system
$$
dq^{i}_{t}=(\beta_{i}-q^{i})\,dt,\quad i=1,...,d.
$$

This fact obviously follows from
Theorems \ref{theorem 4.17.1}
and \ref{theorem 3.8.1} applied to
(\ref{4.24.3})
written in terms of $z^{n}_{t}:=Q^{n}_{0}/n$:
$$
dz^{ni}_{t}=b^{ni}(z^{n}_{t})\,dt+dM^{ni}_{t},
$$
with $d\<M^{n}\>^{ij}_{t}=a^{n}_{t}(z^{n}_{t})\,dt$,
$$
b^{ni}(x)= \delta^{i}(x)\lambda_{0}/n
+\lambda_{i}/n-x , \quad
|a^{nij}_{t}(x)|\leq Nn^{-1}(1+|x|),
$$
where the constant $N$ is independent of $x,n,t$.
}\end{remark}

The following observation can be
generalized so as to  be
used in various control problems
in which optimal controls are discontinuous with respect
to space variables.

\begin{remark}{\rm
                                      \label{remark 4.25.1}
It turns out that many discontinuous functionals
of $x^{n}_{\cdot}$ converge in law to
corresponding functionals of $x_{\cdot}$.
For instance take a Borel vector-valued
function $f(x)$
on $\bR^{d}$ such that the set of its discontinuities
lies in a closed set $J\subset \bR^{d}$ having
Lebesgue measure zero. Also assume that $f$
is locally bounded,
that is bounded on any ball in $\bR^{d}$ but may behave
in any way at infinity.
As an example, one can take   $f(x)=(\delta^{1}(x),
...,\delta^{d}(x))$.
Then, for
$$
y^{n}_{t}:=\int_{0}^{t}f(x^{n}_{s})\,ds,\quad
y_{t}:=\int_{0}^{t}f(x _{s})\,ds
$$
we have that the distributions of $ (x^{n}_{\cdot}
,y^{n}_{\cdot})$ converge weakly to the distribution of
$ (x _{\cdot} ,y _{\cdot})$.

Indeed, append (\ref{4.25.2}) with one more equation:
$dy^{n}_{t}=f(x^{n}_{t})\,dt$ and consider the couple
$z^{n}_{\cdot}=(x^{n}_{\cdot}
,y^{n}_{\cdot})$ as a process in $\bR^{d+1}$.
Obviously Assumptions  \ref{assumption 3.8.2}
and Assumptions  \ref{assumption 3.8.5}
are satisfied for thus obtained couple.

Furthermore, define

$$
H=\{(t,x,y):t>0,y\in\bR,\quad x\in J\quad\text{or}\quad
\prod_{i,j=1}^{d}(\alpha_{i}x^{i}
-\alpha_{j}x^{j})(x^{i}-\nu_{i})=0\}.
$$
Since $J$ is closed, for any $t>0$ and $(x,y)\not\in H_{t}$,
the function $f$ (independent of $y$)
 is continuous in a neighborhood of $x$,
which along with the argument in the proof of Theorem
\ref{theorem 4.18.1} shows that
 Assumption \ref{assumption 3.8.3}
is satisfied for $z^{n}_{t}$. Finally,
for
$$
H^{m}\equiv H,\quad d_{m}=d,\quad v^{mi}(t,x,y)=x^{i},
\quad i=1,...,d,
$$
we have

$$
v^{m}(H_{t})=\{x: x\in J\quad\text{or}\quad
\prod_{i,j=1}^{d}(\alpha_{i}x^{i}
-\alpha_{j}x^{j})(x^{i}-\nu_{i})=0\}
$$
which has $d$-dimensional Lebesgue measure zero and
$$
\det\,V^{nm}(t,x^{n}_{\cdot},y^{n}_{\cdot})=
\det\,a^{n}(x^{n}_{t})\geq
\alpha_{1}^{-1}\cdot...\cdot\alpha_{d}^{-1}>0.
$$
  Hence Assumption \ref{assumption 4.25.2}
is satisfied as well.
This along with precompactness of distributions
of $(x^{n}_{\cdot},y^{n}_{\cdot})$ guaranteed
by Theorem \ref{theorem 4.17.1} and along with
 Theorem \ref{theorem 4.25.1} shows that any convergent
subsequence of distributions
of $(x^{n}_{\cdot},y^{n}_{\cdot})$ converges
to the distribution of a process $(x_{\cdot},y_{\cdot})$,
whose first component satisfies (\ref{4.24.6})
and the second one
obeys $dy_{t}=f(x_{t})\,dt$.

Thus, we get our assertion for a subsequence instead
of the whole sequence.
However, as we have noticed above, solutions
of (\ref{4.24.6}) are weakly unique and this obviously
implies that solutions
of the system (\ref{4.24.6}) appended with
 $dy_{t}=f(x_{t})\,dt$ are also weakly unique.
Therefore, the whole sequence of distributions of
$(x^{n}_{\cdot},y^{n}_{\cdot})$ converges.

}\end{remark}

\section{An $L_{p}$ estimate}
                                     \label{section 3.14.1}

Let $d\geq1$ be an integer,
 $(\Omega,\cF,P)$ be a complete probability space,
and $(\cF_{t},t\geq0)$ be an increasing filtration of
$\sigma$-fields $\cF_{t}\subset\cF$ with $\cF_{0}$
being complete
 with respect to $P,\cF$. Let $K(r,t)$ and $L(r,t)$ be two
nonnegative deterministic function defined for $r,t>0$. Assume
that they increase in $r$ and are locally integrable
in $t$, so that
$$
\int_{0}^{T}(K(r,t)+L(r,t))\,dt<\infty
\quad\forall r,T\in(0,\infty).
$$
Let $\delta(t,x)$ be a nonnegative deterministic function
defined for $t\geq0$ and $x\in\bR^{d}$ and satisfying
$\delta(t,x)\leq K(|x|,t)$.
Define $A(t,x)$
as the set of all  symmetric nonnegative
$d\times d$-matrices
$a$
such that
$$
 \delta(t,x)|\lambda|^{2}\leq
a^{ij}\lambda^{i}\lambda^{j}
\leq K(|x|,t)|\lambda|^{2}\quad\forall\lambda\in\bR^{d}.
$$
Here, as well as  everywhere in the article apart from
 Section \ref{section 4.18.1}, we use the
summation convention.
For any symmetric $d\times d$-matrix $v=(v_{ij})$
define
$$
F(t,x,v)=\sup_{a\in A(t,x)}a^{ij}v_{ij}.
$$
As is easy to see, if $\lambda_{i}(v)$, $i=1,...,d$,
are eigenvalues of $v$ numbered in any order, then
$$
F(t,x,v)=\sum_{i=1}^{d}\chi(t,x,\lambda_{i}(v)),
$$
where $\chi(t,x,\lambda)=K(|x|,t)\lambda$ for $\lambda\geq0$
and $\chi(t,x,\lambda)
=\delta(t,x)\lambda$ for $\lambda\leq0$.

Remember that $C^{\infty}_{0}(\bR^{d+1})$ is
the set of all
infinitely differentiable real-valued function $u=u(t,x)$
on $\bR^{d+1}$ with compact support.

\begin{theorem}
                                     \label{theorem 3.14.1}
Let $x_{t}$, $t\geq0$, be an $\bR^{d}$-valued
$\cF_{t}$-adapted continuous process such that, for any
$u\in C^{\infty}_{0}(\bR^{d+1})$, the following process
 is a local
$\cF_{t}$-supermartingale:
\begin{equation}
                                           \label{2.26.1}
u(t,x_{t})-\int_{0}^{t}\big[
u_{s}(s,x_{s})+F(s,x _{s},u_{x x }(s,x_{s}))
+ L(x^{*}_{t},s)|u_{x}(s,x_{s})|\big]\,ds,
\end{equation}
where $u_{x}$ is the gradient of $u$ with respect to $x$,
$u_{xx}$ is the matrix of second-order derivatives
$u_{x^{i}x^{j}}$ of $u$,
$$
u_{s}=\partial u/\partial s,\quad u_{x^{i}x^{j}}=
\partial^{2}u/\partial x^{i}\partial x^{j}.
$$
 Then for any
$r,T\in(0,\infty)$ there exists a constant $N<\infty$, depending
only on $r,L(r,T)$, and $d$ (but not on $K(r,t)$), such that,
for any nonnegative Borel $f(t,x)$, we have
\begin{equation}
                                             \label{2.26.2}
E\int_{0}^{T\wedge\tau_{r}}\delta^{d/(d+1)}
(t,x_{t})f(t,x_{t})\,dt
\leq N||f||_{L_{d+1}([0,T]\times B_{r})},
\end{equation}
where
$$
||f||_{L_{d+1}([0,T]\times B_{r})}
=\big(\int_{0}^{T}\int_{|x|\leq r}
f^{d+1}(t,x)\,dxdt\big)^{1/(d+1)}.
$$
$B_{r}$ is the open ball in $\bR^{d}$
of radius $r$ centered at the
origin, and $\tau_{r}$ is the first exit time of $x_{t}$
from  $B_{r}$.
\end{theorem}

Proof. First of all notice that
for any $u\in C^{\infty}_{0}(\bR^{d+1})$ expression
(\ref{2.26.1}) makes sense. Indeed, if $r$ is such that
$u(t,x)=0$ for $|x|\geq r$ and all $t$, then the integrand
is bounded by a constant times
$$
\int_{0}^{t}[1+K(r,s)+L(x^{*}_{t}+r,s)]\,ds,
$$
which is finite since
 each trajectory
of $x_{s}$ is bounded on $[0,t]$.
Also observe that usual approximation techniques
allows us to only concentrate on the case
of infinitely differentiable functions $f\geq0$
vanishing for $|x|\geq r$ for some $r$.
We fix $r$, such a function $f$, and a nonnegative function
$\zeta\in C^{\infty}_{0}(\bR^{d+1})$ with unit integral
and support in the  unit ball of $\bR^{d+1}$
centered at the origin.
Below, for any locally bounded Borel function
$g(t,y)$ and $\varepsilon>0$ we use the notation
$$
g^{(\varepsilon)}=g*\zeta_{\varepsilon},\quad
\text{where}\quad\zeta_{\varepsilon}(t,x)=
\varepsilon^{-d-1}\zeta(t/\varepsilon,x/\varepsilon).
$$

Next, we need  Theorem~2 of \cite{Kr76}, which states the
following. There exist constants
$\alpha=\alpha(d)>0$ and
  $N_{r}=N(r,d)<\infty$ and there
exists a bounded Borel nonpositive  function
$z$ on $\bR^{d+1}$
which is convex on $B_{2r} $ for each fixed $t$
and is such that, for  each
nonnegative  symmetric matrix~$a$,
\begin{equation}
                                               \label{2.26.3}
\alpha  (\det a )^{1/(d+1)}f^{(\varepsilon)} \le
  z^{(\varepsilon)}_{t} +
a^{ij} z^{(\varepsilon)}_{x^{i}x^{j}}
\quad\text{for}\quad\varepsilon\leq r,t\in\bR,|x|\leq r,
\end{equation}
\begin{equation}
                                              \label{2.26.5}
|z^{(\varepsilon)}_{x}|\leq2r^{-1}|z^{(\varepsilon)}|
\quad\text{for}\quad\varepsilon\leq r/2,t\in\bR,|x|\leq r,
\end{equation}
\begin{equation}
                                            \label{2.26.4}
 |z|\le N_{r}||f||_{L_{d+1}(\bR\times B_{r})}
\quad\text{in}\quad\bR\times B_{2r}  .
\end{equation}

  Notice that in Theorem~2 of \cite{Kr76}  there  is the minus sign
in front of $z^{(\varepsilon)}_{t}$. However,
(\ref{2.26.3}) is true as is, since one can replace
$t$ with $-t$ and this does not affect any  other term.
Observe that (\ref{2.26.4}) obviously implies that
 for   $\varepsilon\leq r$, we have
\begin{equation}
                                              \label{2.26.6}
 |z^{(\varepsilon)}|\le N_{r}||f||_{L_{d+1}(\bR\times B_{r})}
\quad\text{in}\quad\bR\times B_{r}  .
\end{equation}

Fix an $\varepsilon>0$.
We claim that
the process
\begin{eqnarray}
                                              \label{2.27.1}
\xi_{t}:=-z^{(\varepsilon)}(t\wedge\tau_{r},
x_{t\wedge\tau_{r}})
-\int_{0}^{t\wedge\tau_{r}}\big[-
z^{(\varepsilon)}_{s}(s,x_{s})
\\
+
F(s,x_{s},-z^{(\varepsilon)}_{x x }(s,x_{s}))
+L(r,s)|z^{(\varepsilon)}_{x}(s,x_{s})|\big]\,ds\nonumber
\end{eqnarray}
is a local  supermartingale.
To prove the claim it suffices to prove that
(\ref{2.27.1}) is a local supermartingale on $[0,T]$
for every $T\in[0,\infty)$.
Fix a  $T\in[0,\infty)$
and concentrate on
$t\in[0,T]$. Change $-z^{(\varepsilon)}$ outside of
$[0,T]\times B_{r}$ in any way with the only
requirement that
the  new
function, say $u$ belong to $C^{\infty}_{0}(\bR^{d+1})$.
Then the process (\ref{2.26.1})
is a local supermartingale. Replacing $t$ with
 $t\wedge\tau_{r}$
yields a local supermartingale again. Also observe that
subtracting an increasing
continuous process from a local supermartingale
preserves the property of being a local supermartingale.
After noticing that for $0<s\leq t
\wedge\tau_{r}\leq T$,
we have $|x_{s}|\leq r$ and $L( x^{*}_{s},s) \leq L(r,s)$ and
 we conclude that
$$
\eta_{t}:=u(t\wedge\tau_{r},x_{t\wedge\tau_{r}})
-\int_{0}^{t\wedge\tau_{r}}\big[-
z^{(\varepsilon)}_{s}(s,x_{s})
$$
$$
+F(s,x_{s},-z^{(\varepsilon)}_{xx}(s,x_{s}))
+L(r,s)|z^{(\varepsilon)}_{x}(s,x_{s})|\big]\,ds
$$
is a local supermartingale on $[0,T]$. Since
$$
\eta_{t}-\xi_{t}=[u(0,x_{0})-z^{(\varepsilon)}(0,x_{0})]
I_{\tau_{r}=0},
$$
  is a bounded martingale,
(\ref{2.27.1}) is a local supermartingale indeed.

After having proved our claim
we notice that for each $T\in[0,\infty)$ the process
(\ref{2.27.1}) is obviously bounded on $[0,T]$.
Therefore (\ref{2.27.1}) is a supermartingale and
$$
E\xi_{T}I_{\tau_{r}>0}\leq E\xi_{0}I_{\tau_{r}>0}
\leq\sup_{|x|\leq r}|z^{(\varepsilon)}(0,x)|,
$$
which along with (\ref{2.26.6}), (\ref{2.26.5}),
and the fact that $z\leq0$,
 yields that for any $\varepsilon\leq r/2$
$$
 E\int_{0}^{T\wedge\tau_{r}}\big[
z^{(\varepsilon)}_{s}(s,x_{s})-
F(s,x_{s},-z^{(\varepsilon)}_{xx}(s,x_{s}))\big]\,ds
$$
$$
\leq N_{r}||f||_{L_{d+1}([0,T]\times B_{r})}\big(1
+2r^{-1} E\int_{0}^{T\wedge\tau_{r}}L(r,s)\,ds\big).
$$
Here, owing to (\ref{2.26.3}),
$$
z^{(\varepsilon)}_{s} -
F(s,x,-z^{(\varepsilon)}_{xx} )
=\inf_{a\in A(s,x)}\big[z^{(\varepsilon)}_{s}
+a^{ij}z^{(\varepsilon)}_{x^{i}x^{j}}\big]
$$
$$
\geq f^{(\varepsilon)}\alpha\inf_{a\in A(s,x)}
(\det a)^{1/(d+1)}
=f^{(\varepsilon)}\alpha \delta^{d/(d+1)} .
$$
Hence
$$
E\int_{0}^{t\wedge\tau_{r}}\delta^{d/(d+1)}
f^{(\varepsilon)}(s,x_{s})  \,ds
\leq N ||f||_{L_{d+1}([0,T]\times B_{r})}
$$
with
$$
N=N_{r}\alpha^{-1}\big(1 +2r^{-1}\int_{0}^{T}L(r,s)\,ds\big).
$$

Finally we let $\varepsilon\downarrow0$ and use the continuity
of $f$ which guarantees that $f^{(\varepsilon)}\to f$.
 Then upon remembering that
$f\geq0$ and using Fatou's theorem, we arrive at
(\ref{2.26.2}) with the above specified~$N$.
The theorem is proved.

\begin{remark}{\rm Actually, we did not use the
continuity of $x_{t}$. We could have only assumed that
$x_{t}$ is a separable measurable process. However, then
it turns out that the assumption about the processes
(\ref{2.26.1}) implies that $x_{t}$ is continuous
anyway and  moreover
that $x_{t}$ is an It\^o process (see \cite{Kr1}).
}\end{remark}

\end{document}